\documentstyle[amssymb,amsfonts]{amsart}

\newcommand{\de}{\delta}
\newcommand{\les}{\lesssim}

\newcommand{\nn}{\nonumber}

\renewcommand{\P}{\mbox{{\bf P}}}

\def\bi{\begin{itemize}}
\def\bs{\begin{split}}
\def\es{\end{split}}
\def\ba{\begin{align}}
\def\ea{\end{align}}
\def\eas{\end{align*}}
\def\R{{\hbox{\bf R}}}

\def\Im{{\hbox{\rm Im}}}

\def\diag{{\hbox{diag}}}
\def\R{{\hbox{\bf R}}}
\def\C{{\hbox{\bf C}}}
\def\div{{\hbox{div}}}

\def\Z{{\hbox{\bf Z}}}

\def\eps{\varepsilon}

\newenvironment{proof}{\noindent {\bf Proof} }{\endprf\par}
\def \endprf{\hfill  {\vrule height6pt width6pt depth0pt}\medskip}
\def\emph#1{{\it #1}}
\def\textbf#1{{\bf #1}}
\def\divider#1{$\bullet\quad${\bf #1}}

\def\c{\cdot}

\def\si{\sigma}

\def\f12{\frac{1}{2}}

\def\divider#1{$\bullet\quad${\bf #1}}
 
\def\n{6}		

\theoremstyle{plain}
  \newtheorem{theorem}[subsection]{Theorem}
  
  \newtheorem{proposition}[subsection]{Proposition}
  \newtheorem{lemma}[subsection]{Lemma}

\theoremstyle{remark}

\theoremstyle{definition}
  \newtheorem{definition}[subsection]{Definition}

\include{psfig}

\begin{document}

\title[Global regularity for Maxwell-Klein-Gordon]{Global regularity for the Maxwell-Klein-Gordon equation with small critical Sobolev norm in high dimensions} 
\author{Igor Rodnianski}
\address{Department of Mathematics, Princeton University, Princeton NJ 08544}
\email{ irod@@math.princeton.edu}

\author{Terence Tao}
\address{Department of Mathematics, UCLA, Los Angeles CA 90095-1555}
\email{ tao@@math.ucla.edu}

\subjclass{35Q60\newline
I.R. is a Clay Prize Fellow and supported in part by the NSF grant DMS-01007791.\newline
T.T. is a Clay Prize Fellow and supported in part by a grant from the Packard Foundation}

\vspace{-0.3in}
\begin{abstract}
We show that in dimensions $n \geq \n$ that one has global regularity for the Maxwell-Klein-Gordon equations in the Coulomb gauge provided that the critical Sobolev norm $\dot H^{n/2-1} \times \dot H^{n/2-2}$ of the initial data is sufficiently small.  These results are analogous to those recently obtained for the high-dimensional wave map equation \cite{tao:wavemap1}, \cite{kr:wavemap}, \cite{shsw:wavemap}, \cite{nahmod} but unlike the wave map equation, the Coulomb gauge non-linearity cannot be iterated away directly.  We shall use a different approach, proving Strichartz estimates for the covariant wave equation.  This in turn will be achieved by use of Littlewood-Paley multipliers, and a global parametrix for the covariant wave equation constructed using a truncated, microlocalized Cronstrom gauge.
\end{abstract}

\maketitle

\section{Introduction}

In this paper $n \geq \n$ is an integer, and all implicit constants may depend on $n$.  We let $\R^{1+n}$ be Minkowski space endowed with the usual metric $\eta := \diag(-1,1,\ldots,1)$.

Let $\phi: \R^{1+n} \to \C$ be a complex-valued field, and let $A_\alpha: \R^{1+n} \to \R$ be a one-form; we use $\alpha$ to denote the indices of Minkowski space, which are raised, lowered, and summed in the usual manner.  One can then think of $A$ as a $U(1)$ connection, and can define the \emph{covariant derivatives} $D_\alpha$ by
$$ D_\alpha \phi := (\partial_\alpha + i A_\alpha) \phi.$$
This then induces a \emph{covariant D'Alambertian}
$$ \Box_A := D_\alpha D^\alpha = D^\alpha D_\alpha.$$
We may expand this as
\begin{equation}\label{cov-box}
\Box_A \phi = \Box \phi + 2i A_\alpha \partial^\alpha \phi + i(\partial^\alpha A_\alpha) \phi - (A^\alpha A_\alpha) \phi
\end{equation}
where 
$$ \Box := \partial_\alpha \partial^\alpha = -\partial_t^2 + \Delta$$
is the ordinary D'Alambertian.

We can define the \emph{curvature} $F_{\alpha\beta}$ of the connection $A$ by
$$ F_{\alpha\beta} := \frac{1}{i} [D_\alpha, D_\beta] = \partial_\alpha A_\beta - \partial_\beta A_\alpha.$$

The \emph{Maxwell-Klein-Gordon equations} for a complex field $\phi$ and a one-form $A_\alpha$ are given by
\begin{align*}
\partial^\beta F_{\alpha \beta} &= Im( \phi \overline{D_\alpha \phi} ) \\
\Box_A \phi &= 0;
\end{align*}
these are the Euler-Lagrange equations for the Lagrangian
$$
\int\int  \frac{1}{2} D_{\alpha}\phi \overline{
D^{\alpha} \phi} + \frac{1}{4}F_{\alpha \beta} \overline{F^{\alpha
\beta}}\ dx dt.
$$

The Maxwell-Klein-Gordon system of equations has the gauge invariance
$$
\phi \mapsto e^{i\chi} \phi; \quad A_\alpha \mapsto A_\alpha - \partial_\alpha \chi
$$
for any potential function $\chi: \R^{1+n} \to \R$.  Because of this, it is possible (using Hodge theory) to place this system of equations in the \emph{Coulomb gauge}
\begin{equation}\label{CG}
 \div {\underline{A}} := \partial^j A_j = 0,
\end{equation}
where we use $\underline{A} = (A_j)_{j=1}^n$ to denote the spatial components of $A$, and we use Roman indices to denote summation over the spatial indices $1, \ldots, n$; we will always assume that $\phi$ and $A$ has some decay at spatial infinity, so that this gauge becomes uniquely determined.   In this gauge the Maxwell-Klein-Gordon equations become the following overdetermined system of equations (see e.g. \cite{kl-mac}, \cite{selberg:mkg}, \cite{mac}):
\begin{align*}
\Delta A_0 & = - \Im(\phi \overline{D_0 \phi}) \\
\partial_t \partial_j A_0 &= - (1-\P) \Im (\phi \overline{D_j \phi}) \\
\Box \underline{A}_j &= - \P \Im (\phi \overline{D_j \phi}) \\
\Box'_{\underline{A}} \phi &= 2 iA_0 \partial_t \phi + i(\partial_t A_0) \phi + |\underline{A}|^2 \phi - A_0^2 \phi\\
\div {\underline{A}} &= 0,
\end{align*}
where $\P$ is the Leray projection onto divergence-free vector fields
$$ \P A_k := \Delta^{-1} \partial^j (\partial_j A_k - \partial_k A_j),$$
and $\Box'_{\underline{A}}$ is the modified\footnote{This is of course not a very geometric operator, as we have omitted the terms involving $A_0$ and $\div \underline{A}$, as well as cubic terms of the form $A^2 \phi$.  However, the cubic terms and $A_0$ terms turn out to be negligible, and the $\div \underline{A}$ term vanishes thanks to the Coulomb gauge condition \eqref{CG}.  As it turns out, we will not use much covariant structure in our arguments, although the presence of the $i$ in the lower-order term $2i \underline{A} \cdot \nabla_x$ is crucial in order for a certain phase correction to be purely imaginary.  We should remark however that in three dimensions, recent work of Machedon and Sterbenz \cite{mac} has shown that the full covariant structure of the equation, including the ``elliptic'' component $A_0$ of the gauge, is essential to obtain optimal regularity results.} covariant D'Alambertian
$$ \Box'_{\underline{A}} := \Box + 2i \underline{A} \cdot \nabla_x.$$

We refer to the above system as (MKG-CG).  We observe the well-known fact that the non-linearity for $\Box \underline{A}_j$ has the null form structure\footnote{The term $2i \underline{A} \cdot \nabla_x$ in $\Box'_{\underline{A}}$ also has a null structure due to the Coulomb gauge \eqref{CG}, and this will be exploited to keep various error terms arising from the parametrices we construct for $\Box'_{\underline{A}}$ under control.}
$$ - \P \Im (\phi \overline{D_j \phi}) 
= i\Delta^{-1} \partial_k (\partial_k \phi \overline{\partial_j \phi}
- \partial_j \phi \overline{\partial_k \phi})
+ \P (A_j \phi^2).$$

We shall write $\Phi := (A,\phi)$ to denote the entire collection of fields in (MKG-CG).
Ignoring the tensor structure, constants and Riesz transforms (including the Leray projection $\P$), we can then write (MKG-CG) in caricature form as
\begin{equation}\label{caricature}
\begin{split}
\nabla_x \nabla_{x,t} A_0 & = \phi \nabla_{x,t} \phi + \Phi^3 \\
\Box \underline{A} &= |\nabla_x|^{-1}(\nabla_x \phi \nabla_x \phi) + \Phi^3\\
\Box'_{\underline{A}} \phi &= A_0 \partial_t \phi + (\partial_t A_0) \phi + \Phi^3\\
\div \underline{A} &= 0.
\end{split}
\end{equation}

{\it Remark.} Very crudely, if we ignore the ``elliptic'' component $A_0$, then this equation is roughly of the form
$$ \Box \Phi = \Phi \partial \Phi + \Phi^3.$$
The quadratic non-linearities $\Phi \partial \Phi$ are by far the most difficult to handle, due mainly to the presence of the derivative; the cubic non-linearities $\Phi^3$ are quite manageable and can be handled by the standard theory (based on Strichartz estimates) of semilinear wave equations.

We can study the Cauchy problem for (MKG-CG) by specifying the initial data\footnote{Here and in the sequel we use $\phi[t]$ as short-hand for $(\phi(t), \phi_t(t))$; i.e. $\phi[t]$ encodes both the position and velocity at time $t$.} $\Phi[0]$.  Although we specify initial data for $A_0$, it is essentially redundant since we must obey the compatibility conditions
$$ \div \underline{A}[0] = 0; \quad \Delta A_0(0) = - \Im(\phi(0) \overline{D_0 \phi(0)}); \quad  \partial_j \partial_t A_0 = - (1-\P) \Im (\phi \overline{D_j \phi}).$$

From the scale invariance
$$ \phi(t,x) \mapsto \frac{1}{\lambda} \phi(\frac{t}{\lambda}, \frac{x}{\lambda}); \quad A_\alpha(t,x) \mapsto \frac{1}{\lambda} A_\alpha(\frac{t}{\lambda}, \frac{x}{\lambda})$$
we see that the natural scale-invariant space for the initial data is
$\dot H^{s_c} \times \dot H^{s_c-1}$, where $s_c := n/2-1$ is the critical regularity. 

The Cauchy problem for (MKG-CG) in subcritical regularities $H^s$, $s > s_c$ has been extensively studied (\cite{kl-mac}, \cite{selberg:mkg}, \cite{mac}); we now pause to briefly summarize the known results.  For $s > s_c + 1$ one can obtain local well-posedness by energy methods (see e.g. \cite{em}).  
When $n=3$, this was lowered to $s \geq 1$ in \cite{kl-mac}, and one in fact has global well-posedness in this case by Hamiltonian conservation.  This was improved further to local well-posedness
in $s > 3/4$ in \cite{cuccagna}, and then recently to $s > s_c = \frac{1}{2}$ in \cite{mac}; the global
well-posedness theory was also lowered, but only to $s > \frac{7}{8}$ \cite{keel:mkg}.  In higher dimensions,
local well-posedness was also obtained for $s > s_c$ in \cite{selberg:mkg}; for a simplified model equation for
Maxwell-Klein-Gordon (and Yang-Mills) this is in \cite{kl-tar}.  

At the critical regularity $s = s_c$, much less is known.  When the initial data is small in a certain Besov space\footnote{In the notation of this paper, this result requires $\| \Phi[0] \|_{\dot B^{[2,n]}_1 \times \dot B^{[2,n/2]}_1}$ to be small.}, Sterbenz \cite{sterbenz} 
have recently obtained global well-posedness.  The main result of this paper to replace the Besov norm with the Sobolev norm, at least in high dimensions, and at least for the purposes of obtaining global regularity (rather than well-posedness). 

\begin{theorem}\label{main} If $n \geq \n$ and the initial data $\Phi[0]$ to (MKG-CG) is in $H^s \times H^{s-1}$ for some $s > s_c$, obeys the compatibility conditions, and has a sufficiently small $\dot H^{s_c} \times \dot H^{s_c}$ norm, then the solution $\Phi$ stays in $H^s \times H^{s-1}$ for all time.
\end{theorem}

This is a critical Sobolev regularity result of a flavor similar to the results recently obtained for the wave map equation \cite{tao:wavemap1}, \cite{kr:wavemap}, \cite{nahmod}, \cite{shsw:wavemap}, \cite{tataru}.  As in the theory of the wave map equation, the presence of derivatives in the non-linearity prevents one from obtaining Theorem \ref{main} by a direct iteration argument (i.e. by treating (MKG-CG) as a perturbation of the free wave equation, and iterating away the error in spaces such as Strichartz spaces or $X^{s,b}$ spaces).  In the wave maps equation, this difficulty is resolved by passing to a ``Coulomb gauge'' formulation where the quadratic non-linearity is essentially eliminated, leaving only cubic and better terms.  However this approach is only partly successful in the Maxwell-Klein-Gordon system; the Coulomb gauge \eqref{CG}, while improving the structure of the non-linearity, does not fully renormalize it into a form which can be handled by standard iteration arguments\footnote{In particular, it only generates the ``good'' non-linearity $\nabla_x^{-1}(\nabla_x \phi \nabla_x \phi)$ in the equation for $\underline{A}$.  The equation for $\phi$ still contains bad terms such as $\underline{A} \cdot \nabla_x \phi$ if one fully expands out $\Box_{\underline{A}}$, and this term cannot be iterated away at the critical regularity despite having some null structure, because there is nothing counteracting the possibility of the derivative from landing on a high frequency.  However, the term $\underline{A} \cdot \nabla_x \phi$ is the \emph{only} difficult term in dimensions 5 and higher; if one simply erased this term from (MKG-CG) then one could obtain the global regularity result easily from standard Strichartz estimates.  Indeed, this observation motivates our entire approach, which is to absorb the difficult term into the D'Alambertian and then prove Strichartz estimates for the modified operator. This follows closely the strategy for proving low-regularity existence results for quasilinear wave equations, see e.g. \cite{bc}, \cite{kr:einstein}, \cite{smith} and the references therein.  A similar strategy for the Benjamin-Ono equation was also carried out in \cite{koch}.} based around perturbation of the free d'Alambertian $\Box$. Instead, it appears necessary that one study the (modified) covariant d'Alambertian $\Box'_{\underline{A}}$ directly, and in particular to obtain Strichartz estimates for this operator.  We shall achieve this by using a parametrix which is a modification of the standard Fourier integral operator representation of solution to the free wave equation, but with a phase correction which is inspired by the Cronstrom gauge.

It is quite probable that one can extend this global regularity result to a global small data critical well-posedness result (in analogy with the analogous result for wave maps in \cite{tataru}), and perhaps also obtain scattering as well.  
However, we will not pursue these issues here.  For large data there is unlikely to be any global regularity result due to the energy-supercritical nature of the equation when $n > 4$.  

It would be very interesting to see if this result can be extended to lower dimensions, and especially to the critical case $n=4$ in which the energy is scale invariant.  If one were able to obtain the sharp endpoint Strichartz estimate for covariant wave equations then one could probably also extend to $n=5$ (cf. \cite{tao:wavemap1}, \cite{kr:wavemap};
the point is that in order to get quadratic non-linearities in $L^1_t L^2_x$-type spaces, one needs
the $L^2_t L^4_x$ Strichartz estimate, which is the endpoint in the $n=5$ case).  To extend to $n=4$ however it is likely that one will have to exploit more fully the null structure inherent in the Coulomb gauge \eqref{CG}, perhaps developing covariant bilinear estimates or similar devices.  It is also quite likely that these results extend to Yang-Mills equations; the main differences are that the Coulomb gauge becomes non-linear, and also the connection $A$ does not evolve according to the (perturbed) free wave equation but itself satisfies a covariant wave equation.  

{\bf Acknowledgments}: \,\, The authors are deeply indebted to Sergiu Klainerman, without whose encouragement and insight
this project would not even have been initiated.  We have tremendously benefited from numerous discussions with him, 
in particular on the issue of the Cronstrom gauge. We would also like to thank Joachim 
Krieger and the anonymous referee for their helpful comments and suggestions.

\section{Notation}\label{notation-sec}

We shall need two constants.  First, w need a small number $0 < \delta \ll 1$ depending only on $n$; this measures the loss of exponents in our Strichartz estimates.  Then, we need a sufficiently small number $0 < \eps \ll 1$ depending only on $\delta$ and $n$.  This measures how small our fields $A$, $\phi$ are.

We adopt the convention that $A \lesssim B$ (or $A = O(B)$) denotes the inequality $A \leq CB$ for some constant $C$ depending only on $n$ and $\delta$ (but not on $\eps$). 

If $x = (x_1, \ldots, x_m)$ is a vector, we write $|x| := (\sum_{j=1}^m x_j^2)^{1/2}$ for the $l^2$ norm of $x$.  Similarly for tensors.  We use $\langle x \rangle$ as shorthand for $(1 + |x|^2)^{1/2}$.

We use $\nabla_x$ to denote spatial gradient, and $\nabla_{x,t}$ to denote spacetime gradient.  
We use $|\nabla_x|$ to denote the differentiation operator $|\nabla_x| := \sqrt{-\Delta}$, and $\langle \nabla_x \rangle$ to denote the differentiation operator $\langle \nabla_x \rangle := \sqrt{1-\Delta}$. 

If $\phi(t,x)$ is a function in spacetime, we define the spatial Fourier transform $\hat \phi(t,\xi)$ by
$$ \hat \phi(t,\xi) := \int_{\R^n} e^{-2\pi i x \cdot \xi} \phi(t,x)\ dx.$$
Fix $m(\xi)$ to be a non-negative non-increasing radial bump function supported on $|\xi| \leq 2$ which equals 1 on the ball $|\xi| \leq 1$.  For each integer $k$, we define the Littlewood-Paley projection operators $P_{\leq k} = P_{<k+1}$ to the frequency ball $|\xi| \lesssim 2^k$ by the formula
$$ \widehat {P_{\leq k} \phi}(t,\xi) := m(2^{-k}\xi) \hat \phi(t,\xi),$$
and the projection operators $P_k$ to the frequency annulus $|\xi| \sim 2^k$ by the formula
$$ P_k := P_{\leq k} - P_{<k}.$$
Thus we have $\phi = \sum_k \phi_k$ for test functions $\phi$.  
We also define more general projections $P_{k_1 \leq \cdot \leq k_2}$ by
$$ P_{k_1 \leq \cdot \leq k_2} := P_{\leq k_2} - P_{< k_1}.$$
Similarly define $P_{k_1 < \cdot \leq k_2}$, etc.  

We remark that the Littlewood-Paley  projections defined above commute with all constant-coefficient differential operators and are bounded on every Lebesgue space (including mixed-norm spacetime Lebesgue spaces).  We shall also frequently use the basic product property of Littlewood-Paley operators, in that the expression $P_k(P_{k_1} \phi P_{k_2} \psi)$ vanishes unless we are in one of the following three cases:
\begin{itemize}
\item (High-low interactions) $k_1 = k + O(1)$ and $k_2 \leq k + O(1)$.
\item (Low-high interactions) $k_1 \leq k + O(1)$ and $k_2 = k + O(1)$.
\item (High-high interactions) $k_1 \geq k + O(1)$ and $k_2 = k_1 + O(1)$.
\end{itemize}

We say that a function $\phi$ \hbox{has frequency $\sim 2^k$} if its Fourier transform is supported on the region $|\xi| \sim 2^k$ for all times $t$. Observe that if $\phi$ has frequency $2^k$, then spatial derivatives behave like $2^k$ in the sense that $\| \nabla_x \phi \|_S \sim 2^k \| \phi \|_S$ for any translation invariant Banach space $S$ (basically because the operator $\nabla_x P_k$ has a convolution kernel of mass $O(2^k)$).

We recall \emph{Bernstein's inequality}, which asserts that 
\begin{equation}\label{bernstein}
\| \phi(t) \|_{L^q_x} \lesssim 2^{nk(1/p-1/q)} \| \phi(t) \|_{L^p_x}
\end{equation}
whenever $\phi$ has frequency $2^k$ and $1 \leq p \leq q \leq \infty$.  More generally, we have
$$
\| \phi(t) \|_{L^q_x} \lesssim |V|^{(1/p-1/q)} \| \phi(t) \|_{L^p_x}$$
whenever $\phi$ has Fourier transform supported on a box of volume $V$.

We also need the following standard commutator estimate (cf. \cite{tao:wavemap1}).

\begin{lemma}\label{cook-lemma}  We have
\begin{equation}\label{cook}
\| [P_k, f] g \|_r \lesssim 2^{-k} \| \nabla f\|_{p} \|g\|_{q}
\end{equation}
for integers $k$ and all Schwartz functions $f$, $g$ on $\R^n$ and all $1 \leq p, q, r \leq \infty$ such that $1/p + 1/q = 1/r$.  
\end{lemma}

\begin{proof}
By scaling we may take $k=0$.  We begin with the identity
$$ P_0(fg)(x) - f P_0(g)(x) = \int \check m(y) (f(x-y) - f(x)) g(x-y)\ dy$$
and use the Fundamental theorem of Calculus to rewrite this as
$$- \int_0^1 \int \check m(y) y \cdot \nabla f(x-ty) g(x-y)\ dy\ dt.$$
Since $\check m(y) y$ is integrable, the claim then follows from Minkowski
and H\"older.
\end{proof}

Let $A$ be a connection.  In the following discussion we will assume all functions are sufficiently smooth and decaying at infinity.

For any times $t_1$, $t_0$, we define the linear propagator $S_{\underline{A}}[t_1 \leftarrow t_0]$ by requiring
$$ S_{\underline{A}}[t_1 \leftarrow t_0] \phi[t_0] := \phi[t_1]$$
whenever $\phi$ solves the covariant wave equation $\Box'_{\underline{A}} \phi = 0$.  If $\underline{A}$ is smooth, it is easy to verify (e.g. by classical energy methods) that the Cauchy problem for $\Box'_{\underline{A}}$ is globally well-posed, so that $S_{\underline{A}} [t_1 \leftarrow t_0]$ is well defined.  

For general $\phi$ (not solving the covariant wave equation), we record \emph{Duhamel's formula}
\begin{equation}\label{duhamel-2}
\phi[t] = S_{\underline{A}} [t \leftarrow t_0]\phi[t_0] - \int_{t_0}^t S_{\underline{A}} [t \leftarrow s](0, \Box'_{\underline{A}} \phi(s))\ ds
\end{equation}
for all times $t$, with the convention that $\int_{t_0}^t = - \int_t^{t_0}$ if $t < t_0$.
The Duhamel formula has the consequence that estimates for the homogeneous Cauchy problem can be used to imply similar estimates for the inhomogeneous Cauchy problem.  Indeed, if one has estimates of the form
$$ \| S_{\underline{A}} [t \leftarrow t_0](f,g) \|_X \lesssim \| f \|_Y + \| g \|_Z$$
for all $t_0 \in \R$, where $X, Y, Z$ are Banach spaces defined on $\R \times \R^n$, $\R^n$, and $\R^n$ respectively, then by \eqref{duhamel-2} we see that
$$ \| \phi \|_X \lesssim \| \phi[t_0] \|_{Y \times Z} + \| \Box'_{\underline{A}} \phi \|_{L^1_t Z}.$$

Of course all the above formulae hold for the free wave operator $\Box$, which corresponds to the case $\underline{A}=0$.

\section{Besov spaces}\label{besov-sec}

We now set out our notation for Besov spaces, which will be used heavily in our argument.

We define the inhomogeneous and homogeneous Sobolev spaces $H^s$ and $\dot H^s$ on $\R^n$ by
$$ \| f \|_{H^s} := \| \langle \nabla_x \rangle^s f\|_{L^2(\R^n)}; \quad \| f \|_{\dot H^s} := \| |\nabla_x|^s f \|_{L^2(\R^n)}.$$ 
For any $1 \leq p \leq q \leq \infty$ and $r = 1,2$, we define the homogeneous $l^r$-based Besov space $\dot B^{[p,q]}_r$ by
$$
\| f \|_{\dot B^{[p,q]}_r} := (\sum_{k \in \Z} (2^{(\frac{n}{p} - \frac{n}{q})k} \| P_k f \|_{L^p(\R^n)})^r)^{1/r};$$
the space $\dot B^{[p,q]}_r$ is thus the $L^p$-based homogeneous Besov space which sums $l^r$ in the dyadic frequency pieces, and has the same scaling as $L^q$.  The reason for the notation $[p,q]$ is that Sobolev embedding (or more precisely, Bernstein's inequality) will allow us to control all exponents in the interval $[p,q]$.  Also, the exponent
$q$ will generally be cleaner to write than the regularity $\frac{n}{p}-\frac{n}{q}$, the exact value of which is too
messy to be enlightening in this argument.  When $p=r=2$ we just have $\dot B^{[2,q]}_2 = \dot H^{n/2-n/q}$:
$$ \| f \|_{\dot H^{n/2-n/q}} \sim \| f \|_{\dot B^{[2,q]}_2}.$$
From Bernstein's inequality we have the embedding
\begin{equation}\label{besov-sobolev-scaled} \| f \|_{\dot B^{[\tilde p,q]}_r} \lesssim \| f \|_{\dot B^{[p,q]}_r} 
\end{equation}
whenever $1 \leq p \leq \tilde p \leq q \leq \infty$, and for arbitrary $r$.  Also, from the Littlewood-Paley inequality we have the embedding
\begin{equation}\label{littlewood} \| f \|_{L^p} \lesssim \| f \|_{\dot B^{[p,p]}_2} 
\end{equation}
for all $2 \leq p < \infty$; the condition $p \geq 2$ is necessary in order to interchange the $l^2$ norm of the Littlewood-Paley square function with the $L^p$ norm.

The space $\dot B^{[p,q]}_1$ is a subspace of $\dot B^{[p,q]}_2$, 
\begin{equation}\label{eq:besov}
 \| f \|_{\dot B^{[p,q]}_2} \lesssim \| f \|_{\dot B^{[p,q]}_1},
\end{equation}
and obeys slightly better embedding estimates, in particular $\dot B^{[\infty,\infty]}_1$ embeds into $L^\infty$:
$$ \| f \|_{L^\infty(\R^n)} \lesssim \| f \|_{\dot B^{[\infty,\infty]}_1}.$$
This is an immediate consequence of the triangle inequality.

We observe that the Besov spaces are well-behaved under pseudo-differential operators; indeed for any homogeneous pseudo-differential operator $D^\sigma$ of order $\sigma \in \R$, we have
\begin{equation}\label{ds}
 \| D^\sigma f \|_{\dot B^{[p,\tilde q]}_r} \lesssim \| f \|_{\dot B^{[p,q]}_r}
\end{equation}
whenever $2 \leq p \leq q, \tilde q$ is such that $\frac{1}{\tilde q} = \frac{1}{q} + \frac{\sigma}{n}$.
In particular, zeroth order operators such as the Leray projection $\P$ can be safely ignored if one works exclusively
with Besov spaces.

In particular we have the Sobolev embedding type estimate
\begin{equation}\label{sobolev} \| f \|_{\dot B^{[p,q]}_r} \lesssim \| |\nabla_x|^\sigma f \|_{\dot B^{[p,\tilde q]}_r} 
\end{equation}
under the same assumptions on $p$, $q$, $\tilde q$, $\sigma$.

We now prove a basic product estimate.

\begin{lemma}\label{product}
Whenever $1 \leq p_1 < q_1 < \infty$, $1 \leq p_2 < q_2 < \infty$, and $1 \leq p \le q < \infty$ are such that
$\frac{1}{q} = \frac{1}{q_1} + \frac{1}{q_2}$ and $\frac{1}{p} < \frac{1}{p_1} + \frac{1}{q_2}, \frac{1}{q_1} + \frac{1}{p_2}$, then we have
$$ \| f g \|_{\dot B^{[p,q]}_1} \lesssim \| f \|_{\dot B^{[p_1,q_1]}_2} \| f \|_{\dot B^{[p_2,q_2]}_2}.$$
\end{lemma}

Note that this product estimate improves the $l^2$ Besov norm to an $l^1$ Besov norm; this is ultimately due to the strict inequality in our conditions on $p$.

\begin{proof}
We wish to prove the estimate
$$ \sum_k 2^{(\frac{n}{p}-\frac{n}{q})k} \| P_k(f g) \|_{L^p}
\lesssim
(\sum_{k_1} (2^{(\frac{n}{p_1}-\frac{n}{q_1})k_1} \| P_{k_1} f \|_{L^{p_1}})^2)^{1/2}
(\sum_{k_2} (2^{(\frac{n}{p_2}-\frac{n}{q_2})k_2} \| P_{k_2} g \|_{L^{p_2}})^2)^{1/2}.$$
It will suffice to prove the estimate
$$
2^{(\frac{n}{p}-\frac{n}{q})k} \| P_k(f g) \|_{L^p}
\lesssim
(\sum_{k_1} 2^{(\frac{n}{p_1}-\frac{n}{q_1})k_1} 2^{-\eps|k-k_1|} \| P_{k_1} f \|_{L^{p_1}})
(\sum_{k_2} (2^{(\frac{n}{p_2}-\frac{n}{q_2})k_2} 2^{-\eps|k-k_2|} \| P_{k_2} g \|_{L^{p_2}}))$$
for all $k$ and some absolute constant $\eps > 0$, since the claim then follows by summing in $k$, applying Cauchy-Schwarz, followed by Young's inequality (exploiting the summability of $2^{-\eps|j|}$ in $j$).  By scale invariance (exploiting the condition $\frac{1}{q} = \frac{1}{q_1} + \frac{1}{q_2}$), it suffices to do this for $k=0$.  Splitting $fg = \sum_{k_1} \sum_{k_2} P_{k_1} f_1 P_{k_2} f_2$, it thus suffices to show that
$$
\sum_{k_1,k_2} \| P_0((P_{k_1} f) (P_{k_2} g)) \|_{L^p}
\lesssim
(\sum_{k_1} 2^{(\frac{n}{p_1}-\frac{n}{q_1})k_1} 2^{-\eps|k_1|} \| P_{k_1} f \|_{L^{p_1}})
(\sum_{k_2} (2^{(\frac{n}{p_2}-\frac{n}{q_2})k_2} 2^{-\eps|k_2|} \| P_{k_2} g \|_{L^{p_2}})).$$
There are only three cases in which the summand is non-zero: when $k_1 = O(1)$ and $k_2 \leq O(1)$ (the ``high-low'' case); when $k_1 \leq O(1)$ and $k_2 = O(1)$ (the ``low-high'' case); and when $k_1 = k_2 + O(1) \geq O(1)$ (the ``high-high'' case).

Consider first the contribution of the high-low case.  In this case we discard the bounded $P_0$ multiplier and use H\"older to estimate
$$ \| P_0((P_{k_1} f) (P_{k_2} g)) \|_{L^p} \lesssim \| P_{k_1} f \|_{p_1} \| P_{k_2} f \|_{r_2}$$
where $1/r_2 := 1/p - 1/p_1$. By hypothesis, we have $r_2 > q_2$.  Hence by Bernstein's inequality \eqref{bernstein},
$$ \| P_{k_2} f \|_{r_2} \lesssim 2^{(\frac{n}{p_2}-\frac{n}{q_2})k_2} 2^{-\eps|k_2|} \| P_{k_2} g \|_{L^{p_2}})$$
for some $\eps > 0$ (here we are exploiting the fact that $k_2 \leq O(1)$).  Since $k_1 = O(1)$, we thus see that the high-low contribution is acceptable.

The contribution of the low-high case can be dealt with similarly, so we now turn to the high-high case.  Observe that in the high-high case it suffices to prove the estimate
$$
\|P_0 ((P_{k_1} f) (P_{k_2} g))\|_{L^p}\lesssim  2^{(\frac n{p_1}-\frac n{q_1})k}
\|P_{k_1} f\|_{L^{p_1}} 2^{(\frac n{p_2}-\frac n{q_2})k} \|P_{k_2} f\|_{L^{p_2}} 
$$
Choose  $p_1 < r_1 \le q_1$, $p_2 < r_2 \le q_2$ such that $1/p = 1/r_1 + 1/r_2$; this is possible thanks to our assumptions on $q,q_1,q_2$.  By Bernstein and H\"older we thus have
$$ \| P_0((P_{k_1} f) (P_{k_2} g)) \|_{L^p} \lesssim 
\| P_{k_1} f \|_{L^{r_1}} \| P_{k_2} f \|_{L^{r_2}}.$$
But by Bernstein again we have
$$ \| P_{k_j} f \|_{L^{r_j}} \lesssim 2^{(\frac{n}{p_j}-\frac{n}{q_j})k_j}  \| P_{k_j} f\|_{L^{p_j}}
$$
for $j=1,2$ (here we are exploiting the fact that $k_j \leq O(1)$ and $r_j \le q_j$).  Thus the high-high contribution is also acceptable.
\end{proof}

By applying Lemma \ref{product} twice (and using \eqref{eq:besov}), we see that we have the trilinear estimate
\begin{equation}\label{eq:trilinear}
 \| f g h \|_{\dot B^{[p,q]}_1} \lesssim \| f \|_{\dot B^{[p_1,q_1]}_2} \| g \|_{\dot B^{[p_2,q_2]}_2} \| h\|_{\dot B^{[p_3,q_3]}_2}
\end{equation}
whenever $1 \leq p_j < q_j < \infty$ for $j=1,2,3$ and $1 \leq p \le q < \infty$ are such that $\frac{1}{q} = \frac{1}{q_1} + \frac{1}{q_2} + \frac{1}{q_3}$ and
$$ \frac{1}{p} < \frac{1}{p_1} + \frac{1}{q_2} + \frac{1}{q_3}, \frac{1}{q_1} + \frac{1}{p_2} + \frac{1}{q_3}, \frac{1}{q_1} + \frac{1}{q_2} + \frac{1}{p_3}.$$

We will also need the following variant of Lemma \ref{product}:

\begin{lemma}\label{prod-variant} We have the spacetime estimate
$$ \| F G \|_{L^1_t \dot B^{[2,n/2]}_2}
\lesssim \|F\|_{L^1_t \dot B^{[\infty,\infty]}_1} \|G\|_{L^\infty_t \dot B^{[2,n/2]}_2}
+ \| F \|_{L^2_t \dot B^{[p,2n]}_2} \| G \|_{L^2_t \dot B^{[q,2n/3]}_2}$$
whenever $2 \leq p <2n$,  $p\le 2n/(n-3)$ and $2 \leq q < 2n/3$. Note that  
the conditions on $p,q$ imply that 
$1/p+1/q > 1/2$. 
\end{lemma}

\begin{proof}
We split $FG = \sum_{k_1} \sum_{k_2} (P_{k_1} F) (P_{k_2} G)$.  It will suffice to prove the paraproduct estimates
\begin{align*}
\| \sum_{k_1, k_2: k_1 \leq k_2 + C} (P_{k_1} F) (P_{k_2} G) \|_{L^1_t \dot B^{[2,n/2]}_2} &\lesssim
\|F\|_{L^1_t \dot B^{[\infty,\infty]}_1} \|G\|_{L^\infty_t \dot B^{[2,n/2]}_2} \\
\| \sum_{k_1, k_2: k_1 > k_2 + C} (P_{k_1} F) (P_{k_2} G) \|_{L^1_t \dot B^{[2,n/2]}_2} &\lesssim
\| F \|_{L^2_t \dot B^{[p,2n]}_2} \| G \|_{L^2_t \dot B^{[q,2n/3]}_2}.
\end{align*}
We begin with the first inequality.  By a H\"older in time it suffices to prove the spatial estimate
$$\| \sum_{k_1, k_2: k_1 \leq k_2 + C} \sum_{k_2} (P_{k_1} f) (P_{k_2} g) \|_{\dot B^{[2,n/2]}_2} \lesssim
\|f\|_{\dot B^{[\infty,\infty]}_1} \|g\|_{\dot B^{[2,n/2]}_2}$$
for any functions $f, g$.  By the definition of the $\dot B^{[\infty,\infty]}_1$ norm and the triangle inequality, it suffices to show
$$\| \sum_{k_2: k_1 \leq k_2 + C}  (P_{k_1} f) (P_{k_2} g) \|_{\dot B^{[2,n/2]}_2} \lesssim
\|P_{k_1} f\|_{\infty} \|g\|_{\dot B^{[2,n/2]}_2}$$
for each $k_1$.  By scale invariance we may take $k_1 = 0$.

If $k_2 > C$ then the functions $(P_{k_1} f) (P_{k_2} g)$ essentially have frequency $2^{k_2}$, and in particular are orthogonal.  The claim then follows by expanding out the $\dot B^{[2,n/2]}_2$ norm and taking $P_{k_1} f$ out in $L^\infty$.  If $-C \leq k_2 \leq C$ then the function $(P_{k_1} f) (P_{k_2} g)$ has Fourier support on the region $|\xi| = O(1)$, and so one can estimate the $\dot B^{[2,n/2]}_2$ norm by the $L^2$ norm.  The claim then follows from H\"older's inequality (since the $L^2$ norm of $P_{k_2} g$ is comparable to its $\dot B^{[2,n/2]}_2$ norm.

Now we prove the second inequality.  By another H\"older in time, it suffices to show that
$$
\| \sum_{k_1, k_2: k_1 > k_2 + C} (P_{k_1} f) (P_{k_2} g) \|_{\dot B^{[2,n/2]}_2} \lesssim
\| f \|_{\dot B^{[p,2n]}_2} \| g \|_{\dot B^{[q,2n/3]}_2}.
$$
Observe that the expression $(P_{k_1} f) (P_{k_2} g)$ has frequency $\sim 2^{k_1}$.  Thus it will suffice to prove that
$$
2^{\frac{n-4}{2} k_1} \| \sum_{k_2: k_1 > k_2 + C} (P_{k_1} f) (P_{k_2} g) \|_{L^2} \lesssim
2^{(\frac{n}{p}-\frac{1}{2})k_1} \| P_{k_1} f \|_{L^p} \| g \|_{\dot B^{[q,2n/3]}_2}
$$
for each $k_1$, since the claim then follows by square-summing in $k_1$.  By scaling again we may take $k_1 = 0$.  But since $1/p + 1/q > 1/2$, we see from Bernstein and H\"older that
$$ \| (P_0 f) (P_{k_2} g) \|_{L^2} \lesssim
\| P_0 f \|_{L^p} \| P_{k_2} g \|_{L^{\frac {2p}{p-2}}}\lesssim 2^{(\frac{3}{2} - 
\frac{n(p-2)}{2p})k_2} 2^{-\epsilon |k_2|}
\| P_0 f \|_{L^p} \| g \|_{\dot B^{[q,2n/3]}_2}.$$
Summing in $k_2\le O(1)$ we obtain the result with the help of the condition that 
$p\le 2n/(n-3)$.
\end{proof}

\section{Iteration spaces}\label{function-sec}

We define some Banach spaces $D$, $S$, $S_{ellip}$, $N_1$, $N_2$ which we will iterate in.  Our initial data $\Phi[0]$ will be measured in the data norm $D$, defined by
\begin{equation}\label{data}
\| \Phi[0] \|_{D} := \| \nabla_{x,t} \Phi(0) \|_{\dot B^{[2,n/2]}_2}
\sim \| \Phi[0] \|_{\dot H^{n/2-1}_x \times \dot H^{n/2-2}_x}.
\end{equation}
Similarly
the solution $\Phi$ will be measured in the solution norm $S$, defined by
\begin{equation}\label{solution}
\| \Phi \|_{S} := \| \nabla_{x,t} \Phi \|_{L^\infty_t \dot B^{[2,n/2]}_2} + 
\| \nabla_{x,t} \Phi \|_{L^2_t \dot B^{[p_*,2n/3]}_2}
\end{equation}
where $p_* := \frac{2(n-1)}{n-3} + \delta$, and $0 < \delta \ll 1$ is the small number depending on $n$ chosen previously.  Note that the hypothesis $n \geq \n$ allows us to choose $\delta$ so that $p_* < 2n/3$.  The significance of $p_*$ is that it is slightly bigger than the endpoint Strichartz exponent.

The elliptic portion $A_0$ of the field $\Phi$ will be measured in the elliptic solution norm $S_{ellip}$, defined by
\begin{equation}\label{elliptic}
\| A_0 \|_{S_{ellip}} := \| \nabla_{x,t} A_0 \|_{L^1_t \dot B^{[p_{**},n]}_1} + \| \nabla_{x,t} A_0 \|_{L^\infty_t \dot B^{[2,n/2]}_2},
\end{equation}
where $p_{**} := \frac{2n}{n-2} - \delta$; note that the first Besov space is summed in $l^1$ rather than $l^2$.
Finally, the non-linearity will be measured either in the norm $N_2$, defined by
\begin{equation}\label{nonlinear}
\| G \|_{N_2} := \| G \|_{L^1_t \dot B^{[2,n/2]}_2} \sim \| G \|_{L^1_t \dot H^{n/2-2}_x},
\end{equation}
or in the space $N_1$, defined by
\begin{equation}\label{nonlinear-1}
\| G \|_{N_1} := \| G \|_{L^1_t \dot B^{[2,n/2]}_1}.
\end{equation}
Clearly the $N_1$ norm controls the $N_2$ norm:
$$ \| G\|_{N_2} \lesssim \| G\|_{N_1}.$$

Let $p_{***}$ be the exponent\footnote{The reader may wish to use the concrete case $n=6$ to track the numerology; in this case $p_* \approx 10/3$, $p_{**} \approx 12/5$, and $p_{***} \approx 3$.} defined by $\frac{1}{p_{***}} =\frac{1}{2} ( \frac{1}{2} + \frac{1}{p_{**}} )$; thus $p_{***}$ is slightly less than $2n/(n-1)$.  From interpolation and Sobolev embedding we observe that 
\begin{equation}\label{eq:psss}
\| \nabla_{x,t} A_0 \|_{L^2_t \dot B^{[p_{***},2n/3]}_2} \lesssim \| A_0 \|_{S_{ellip}},
\end{equation}
which implies in particular that $S_{ellip}$ contains $S$:
\begin{equation}\label{e-s}
\| A_0 \|_S \lesssim \| A_0 \|_{S_{ellip}}.
\end{equation}

The spaces $D$, $S$, and $S_{ellip}$ have the scaling of $length^{-1}$, while $N$ has the scaling of $length^{-3}$.  Thus if $\Phi \in S$, we expect $N$ to contain such quantities as $\Box \Phi$, $\Phi \nabla_x \Phi$ or $\Phi^3$.

Clearly $S$ controls $D$ on time slices:
\begin{equation}\label{s-d}
\sup_{t \in I} \| \Phi[t] \|_D \lesssim \| \Phi \|_S.
\end{equation}
We also have the standard Strichartz estimate\footnote{As we are in the high-dimensional case $n \geq \n$, we will not need the endpoint Strichartz estimate $L^2_t L^{2(n-1)/(n-3)}$ here.  Indeed we will rely on this room in the exponents when we prove covariant Strichartz estimates later in this paper; when $n=5$ it seems that one is forced to
resort to either endpoint covariant Strichartz estimates, or covariant bilinear estimates, to recover the global
regularity results.} 
\begin{equation}\label{strichartz}
\| \Phi \|_S \lesssim \| \Phi[t_0] \|_D + \| \Box \Phi \|_{N_2}
\end{equation}
for any time $t_0$; see e.g. \cite{tao:keel}.

\section{The covariant Strichartz estimate}\label{covariant-sec}

The proof of Theorem \ref{main} will rely crucially on the following (modified-)covariant Strichartz estimate.

\begin{proposition}[$\dot H^{n/2-1}$ covariant Strichartz estimate]\label{covariant-Strichartz}  Let $t_0$ be a time, let $I$ be any compact time interval containing $t_0$, and let $A$ be a smooth connection on $I \times \R^n$ which obeys the smallness condition
$$
\| \underline{A} \|_D + \| \Box \underline{A} \|_{N_1} \lesssim \eps
$$
and the Coulomb gauge condition 
\begin{equation}\label{CG-underline}
\div \underline{A} = 0
\end{equation}
on $I \times \R^n$.  Then we have the (modified) covariant Strichartz estimates
\begin{equation}\label{strichartz-covariant}
\| \phi \|_S \lesssim \| \phi[t_0] \|_D + \| \Box'_{\underline{A}} \phi \|_{N_2}
\end{equation}
on $I \times \R^n$ and all Schwartz $\phi$ on $I \times \R^n$.  (Note that all the implicit constants are independent of $I$).
\end{proposition}

{\it Remark.}  This estimate is clearly a generalization of the standard Strichartz estimate \eqref{strichartz}, which is the special case $\underline A=0$.  Roughly speaking, it asserts that the modified operator $\Box'_{\underline{A}}$ can
be treated as if it were equivalent to $\Box$ for the purposes of proving global regularity via Strichartz estimates.  Note that the time component $A_0$ of the connection is irrelevant here since it does not appear in the modified covariant D'Alambertian $\Box'_{\underline{A}}$.  Also note that we are assuming that $\Box \underline{A}$ is small with respect to the $l^1$-based Besov space $N_1$, and not just the more familiar Sobolev space $N_2$.  This stronger assumption will be important in our argument.  The Coulomb gauge assumption \eqref{CG-underline}
is crucial to our argument in controlling error terms, although it seems that for sufficiently high dimension
(e.g. $n \geq 10$) it is possible to obtain acceptable control on all error terms without requiring the Coulomb
gauge assumption.

For the rest of this section we shall assume Proposition \ref{covariant-Strichartz} and show how it implies Theorem \ref{main}.   Then in the remainder of the paper we shall prove Proposition \ref{covariant-Strichartz}.

We now begin the proof of Theorem \ref{main}.  The idea will be to treat all the terms in the right-hand side of \eqref{caricature} as negligible error terms, using the Strichartz estimates for both $\Box$ and $\Box'_{\underline{A}}$ to do this; the numerology will allow us to do this because we are in high dimensions $n \geq 6$ and because the derivatives in the non-linearities on the right-hand sides of (MKG-CG) are in favorable locations.  Proposition \ref{covariant-Strichartz} in turn allows us to obtain these Strichartz estimates assuming that the right-hand sides were indeed small, thus closing the bootstrap argument. 

By time reversal symmetry it suffices to prove uniform $H^s \times H^{s-1}$ bounds on an arbitrary time interval $I = [0,T]$, which we now fix; our bounds will be independent of $T$.  We remark from the existing local well-posedness theory (see e.g. \cite{selberg:mkg}) that we may assume \emph{a priori} that $\Phi$ is smooth and has some decay at infinity. 

Fix $I$, $\Phi$, $s$.  From hypothesis we may assume that
\begin{equation}\label{hypothesis}
 \| \Phi[0] \|_D \lesssim \eps.
\end{equation}

The key estimate to prove in our argument is the \emph{a priori} estimate
\begin{equation}\label{se-control}
\| \Phi \|_S + \| A_0 \|_{S_{ellip}} \lesssim \eps.
\end{equation}
Note that the right-hand side is independent of the time interval $I$.

We now prove \eqref{se-control}.  By a simple continuity argument we may assume as a bootstrap hypothesis that
\begin{equation}\label{se-control-2}
\| \Phi \|_S + \| A_0 \|_{S_{ellip}} \lesssim \eps^{1/2}
\end{equation}
if $\eps$ is chosen sufficiently small.

Consider first the contribution of $A_0$.  We begin with the $L^\infty_t \dot B^{[2,n/2]}_2$ component of the $S_{ellip}$ norm, which is an easy term.  By \eqref{caricature} it suffices to show that
$$ \| |\nabla_x|^{-1} (\phi \nabla_{x,t} \phi + \Phi^3) \|_{L^\infty_t \dot B^{[2,n/2]}_x}
\lesssim \eps.$$

By \eqref{ds}, the left-hand side can be bounded by
$$ \lesssim \| \phi \nabla_{x,t} \phi + \Phi^3 \|_{L^\infty_t \dot B^{[2,n/3]}_2}.$$
By several applications of Lemma \ref{product}, we may bound this by
$$ \lesssim \| \phi \|_{L^\infty_t \dot B^{[2,n]}_2} \| \phi \|_{L^\infty_t \dot B^{[2,n/2]}_2}
+ \| \Phi \|_{L^\infty_t \dot B^{[2,n]}_2}^3.$$
But by \eqref{sobolev}, \eqref{se-control-2} this is $O( (\eps^{1/2})^2 + (\eps^{1/2})^3) = O(\eps)$, which is acceptable.

Now we consider the $L^1_t \dot B^{[p_{**},n]}_1$ norm, which is also fairly easy.  Arguing as before, we have to show that
$$ \| \phi \nabla_{x,t} \phi + \Phi^3 \|_{L^1_t \dot B^{[p_{**},n/2]}_1} \lesssim \eps.$$
By Lemma \ref{product} and \eqref{eq:trilinear}, we may bound the left-hand side by
$$ \lesssim \| \phi \|_{L^2_t \dot B^{[p_*, 2n]}_2} \| \nabla_{x,t} \phi \|_{L^2_t \dot B^{[p_*, 2n/3]}_2}
+ \| \Phi \|_{L^2_t \dot B^{[p_*, 2n]}_2}^2 \| \Phi \|_{L^\infty_t \dot B^{[2,n]}_2},$$
where we have used the assumption $n \geq \n$ (in a rather weak way) to ensure that $p_*$ is large enough for the hypotheses of Lemma \ref{product} to be respected.
This completes the estimation of $A_0$ in $S_{ellip}$ for \eqref{se-control}.

It remains to estimate the $\|\Phi\|_S$ component of \eqref{se-control}.  We first need a preliminary estimate.

\begin{lemma}\label{boxa-lemma}  We have
$$ \|\Box \underline{A}\|_{N_1}\lesssim \eps.$$
\end{lemma}

\begin{proof} 
By \eqref{caricature} it suffices to prove that
$$
\| |\nabla_x|^{-1} ((\nabla_x \phi)^2)) \|_{N_1}
+ \| \Phi^3 \|_{N_1}
\lesssim \eps.$$
We first deal with the quadratic term $\| |\nabla_x|^{-1} (\nabla_x \phi \nabla_x \phi) \|_{N_1}$.  By \eqref{ds}, we can estimate this by
$$\| |\nabla_x|^{-1} (\nabla_x \phi \nabla_x \phi) \|_{N_1} \lesssim \| (\nabla_x \phi)^2 \|_{L^1_t \dot B^{[2,n/3]}_1};$$
note we definitely need the hypothesis $n \geq \n$ here in order to keep the regularity non-negative.
By Lemma \ref{product}, we thus have
$$ \| |\nabla_x|^{-1} (\nabla_x \phi \nabla_x \phi) \|_{N_1} \lesssim \| \nabla_x \phi \|_{L^2_t \dot B^{[p_*,2n/3]}_2}^2;$$
note that the hypotheses of Lemma \ref{product} will be satisfied since we have $\frac{1}{2} < \frac{1}{p_*} + \frac{3}{2n}$ when $n \geq \n$ (and $\delta$ is chosen sufficiently small).  By \eqref{se-control-2} (and \eqref{sobolev}) this in turn is bounded by $O(\eps)$, as desired.

Now we turn to the cubic term.  By \eqref{eq:trilinear} we will have
$$ \| \Phi^3 \|_{L^1_t B^{[2,n/2]}_1} \lesssim \| \Phi \|_{L^2_t B^{[p_*,2n]}_2}^2 \| \Phi \|_{L^\infty_t \dot B^{[2,n]}_2}$$
provided that we have
$$ \frac{1}{2} < \frac{1}{p_*} + \frac{1}{2n} + \frac{1}{n}, \frac{1}{2} + \frac{1}{2n} + \frac{1}{2n}.$$
But these inequalities certainly hold for $n \geq \n$ (if $\delta$ is sufficiently small).   By \eqref{se-control-2} and \eqref{sobolev} we may thus bound 
$$\| \Phi^3 \|_{N_1} \lesssim (\eps^{1/2})^3 \lesssim \eps$$
as desired.
\end{proof}

We now return to the estimation of $\|\Phi\|_S$.  By \eqref{e-s} we see that the $A_0$ component of $\Phi$ is already satisfactorily estimated, so we focus on $\underline{A}$ and $\phi$.  The estimate for $\underline{A}$ follows from \eqref{strichartz}, \eqref{hypothesis}, and Lemma \ref{boxa-lemma}, so it suffices to control $\phi$.

From \eqref{hypothesis} and Lemma \ref{boxa-lemma} we see that $A$ obeys all the hypotheses required to invoke Proposition \ref{covariant-Strichartz}.  It thus suffices to prove the estimates 
$$
\| A_0 \partial_t \phi \|_{N_2}
+ \| (\partial_t A_0) \phi \|_{N_2}
+ \| \Phi^3 \|_{N_2}
\lesssim \eps.$$
The cubic term was already proven to be acceptable in the proof of Lemma \ref{boxa-lemma}, so we turn to the quadratic terms.  First consider the contribution of $A_0 \partial_t \phi$.  By Lemma \ref{prod-variant} we have
$$ \| A_0 \partial_t \phi \|_{L^1_t \dot B^{[2,n/2]}_2}
\lesssim \| A_0 \|_{L^1_t \dot B^{[\infty,\infty]}_1}
\| \partial_t \phi \|_{L^\infty_t L^{[2,n/2]}_2}
+
\| A_0 \|_{L^2_t \dot B^{[p_{***},2n]}_2}
\| \partial_t \phi \|_{L^2_t L^{[p_*,2n/3]}_2};$$
note that $p_{***}\le \frac{2n}{n-3}$ and $p_*< 2n/3$ when $n \geq \n$ and $\delta$ is sufficiently small.  Applying \eqref{se-control-2}, \eqref{sobolev}, and \eqref{eq:psss} we thus obtain
$$ \| A_0 \partial_t \phi \|_{L^1_t \dot B^{[2,n/2]}_2}
\lesssim \| A_0 \|_{S_{ellip}} \| \phi \|_S \lesssim \eps$$
as desired.

Now we consider the contribution of $(\partial_t A_0) \phi$.  For this we just use Lemma \ref{product} to estimate
$$ \| (\partial_t A_0) \phi \|_{L^1_t \dot B^{[2,n/2]}_2}
\lesssim \| \partial_t A_0 \|_{L^2_t \dot B^{[p_{***},2n/3]}_1}
\| \phi \|_{L^2_t \dot B^{[p_*,2n]}_2};$$
note that the condition $n \geq \n$ implies that
$$ \frac{1}{2} < \frac{1}{p_{***}} + \frac{1}{2n}, \frac{3}{2n} + \frac{1}{p_*}.$$
Applying \eqref{se-control-2}, \eqref{sobolev}, \eqref{eq:psss} we thus obtain
$$ \| (\partial_t A_0) \phi \|_{L^1_t \dot B^{[2,n/2]}_2}
\lesssim \| A_0 \|_{S_{ellip}} \| \phi \|_S \lesssim \eps$$
as desired. This completes the proof of \eqref{se-control}.

We now use \eqref{se-control} to prove the $H^s \times H^{s-1}$ regularity.  Since this type of argument is well known we only provide a sketch of it here\footnote{An alternate way to proceed at this point is to use frequency envelopes, as in \cite{tao:wavemap1}.}.  For simplicity we deal with the case when $s = s_c + 1$, although the other cases are similar (one needs to use fractional Leibnitz instead of integer Leibnitz).  

By assumption we may assume that
$$ \| \nabla_x \Phi[0] \|_D \lesssim M$$
for some finite constant $M$.  By \eqref{e-s} it will suffice to show that
\begin{equation}\label{eq:continuity}
 \| \nabla_x \Phi \|_S + \| \nabla_x A_0 \|_{S_{ellip}} \lesssim M
\end{equation}
for all $0 \leq \alpha \leq n/2-s$.  As before we may use continuity arguments to assume a priori that
\begin{equation}\label{continuity-2}
\| \nabla_x \Phi \|_S + \| \nabla_x A_0 \|_{S_{ellip}} \lesssim \eps^{-1/2} M.
\end{equation}

One now differentiates \eqref{caricature} by $\nabla_x$ and repeats the previous analysis.  From the ordinary Leibnitz rule one of the factors on the right-hand side will acquire a $\nabla_x$, and this term will be estimated using \eqref{eq:continuity} instead of \eqref{se-control-2}; by repeating the previous estimates we will obtain \eqref{eq:continuity}.  The only interesting case occurs is when one commutes $\nabla_x $ with the covariant D'Alambertian
$\Box'_{\underline{A}}$, as one picks up additional terms of the form
$(\nabla_x \underline{A}) (\nabla_x \phi)$.  But these are similar to the terms $\nabla( |\nabla_x|^{-1}((\nabla_x \phi) (\nabla_x \phi)))$ that one will pick up anyway from differentiating \eqref{caricature}, and are estimated in the same way.  We omit the details.

It now remains only to prove the covariant Strichartz estimate in Proposition \ref{covariant-Strichartz}.  This will occupy the remainder of the paper.

\section{Reduction to a frequency-localized parametrix }\label{localized-sec}

We now begin the proof of the covariant Strichartz estimate.  The first step is to reduce matters to the following frequency-localized Strichartz estimate for a parametrix, where the wave $\phi$ is constrained to have higher frequency than the connection $\underline{A}$.  The main reason that we can obtain this reduction is that all the terms of the form $\underline{A} \cdot \nabla_x \phi$ in which the frequency of $\underline{A}$ has higher frequency of $\phi$ are manageable, because the derivative falls on a low frequency term.

\begin{proposition}[Frequency-localized covariant Strichartz estimate for a parametrix]\label{localized-Strichartz}  Let $t_0$ be a time, let $I$ be a compact time interval containing $t_0$, and let $A$ be a smooth connection on $I \times \R^n$ which obeys the smallness condition
$$
\| \underline{A} \|_D + \| \Box \underline{A} \|_{N_1} \lesssim \eps
$$
and the Coulomb gauge condition \eqref{CG-underline}
on $I \times \R^n$.  
Suppose also that $\underline{A}$ has Fourier support on the region $\{ \xi: |\xi| \leq 2^{k-10} \}$ for some integer $k$.  Let $(f,g)$ be a pair of Schwartz functions on $\R^n$, and $F$ be a function on $I \times \R^n$,
all with Fourier support in the region $\{ \xi: 2^{k-3} \leq |\xi| \leq 2^{k+3} \}$, and with norm
$$ \| (f,g) \|_D + \| F \|_{N_2} = K$$
for some $K > 0$.   Then there exists a function $\phi$ on $I \times \R^n$ 
with Fourier support in the region $\{ \xi: 2^{k-10} \leq |\xi| \leq 2^{k+10} \}$ obeying the estimates 
\begin{align*}
\| \phi \|_{S} &\lesssim K,\\
\| \phi[t_0] - (f,g) \|_D &\lesssim \eps^{\delta} K,\\
\| \Box'_{\underline{A}} \phi - F \|_{N_2} &\lesssim \eps^{\delta} K
\end{align*}
on $I \times \R^n$.
\end{proposition}

We now show how the above Proposition implies the covariant Strichartz estimate \eqref{strichartz-covariant}.

Fix $I$, all our computations below shall be restricted to the spacetime slab $I \times \R^n$.  Let $\underline{A}$ be as in Proposition \ref{covariant-Strichartz}.  From \eqref{strichartz} we have in particular that
\begin{equation}\label{eq:a-s}
\| \underline{A} \|_S \lesssim \eps.
\end{equation}
We have to prove the estimate \eqref{strichartz-covariant}.  Since $\underline{A}$ is smooth and Schwartz, and $I$ is compact, it is easy to see (thanks to the ordinary Strichartz estimate \eqref{strichartz} and perturbation theory) that we have some estimate of the form
$$ \| \phi \|_S \lesssim C(\underline{A}, I) (\| \phi[t_0] \|_D + \| \Box'_{\underline{A}} \phi \|_{N_2}),$$
and furthermore the constant $C(\underline{A},I)$ depends continuously on $I$ and on $\underline{A}$ in a smooth topology.  Thus by a continuity argument, to prove \eqref{strichartz-covariant} it will suffice to do so under the a priori assumption that
\begin{equation}\label{strichartz-covariant-bootstrap}
\| \phi \|_S \lesssim \eps^{-\delta/2} (\| \phi[0] \|_D + \| \Box'_{\underline{A}} \phi \|_{N_2})
\end{equation}
for all smooth Schwartz $\phi$ on $I \times \R^n$.

Now we prove \eqref{strichartz-covariant}.  From the remarks following the Duhamel formula \eqref{duhamel-2}, we see that it suffices to prove the estimate
$$ \| S_{\underline{A}}[t \leftarrow t_0] (f,g) \|_S \lesssim \| (f,g) \|_D$$
for all times $t_0 \in I$ and all Schwartz functions $f,g$.

Fix $t_0, f, g$; we may normalize $\| (f,g) \|_D = 1$.  Write $\phi(t) := S_{\underline{A}}[t \leftarrow t_0] (f,g)$; thus $\phi[t_0] = (f,g)$, and $\Box'_{\underline{A}} \phi = 0$.  By \eqref{strichartz-covariant-bootstrap} we have
\begin{equation}\label{eq:phis-big}
 \| \phi \|_S \lesssim \eps^{-\delta/2}.
\end{equation}

Our goal is to improve this to
\begin{equation}\label{phis:targ}
 \| \phi \|_S \lesssim 1
\end{equation}
We use the Littlewood-Paley operators to split $\phi = \sum_k \phi_k$.  By \eqref{solution} and the definition of the $l^2$-based Besov spaces $\dot B^{s,p}_2$, we see that
$$ \| \phi \|_S \lesssim (\sum_k \| P_k \phi \|_S^2)^{1/2}.$$
For each $k$, we apply Proposition \ref{localized-Strichartz} with $\underline{A}$ replaced by $P_{\leq k-20} \underline{A}$, $(f,g)$ set equal to $P_k \phi[t_0]$, and $F$ set equal to $\Box'_{P_{\leq k-20} \underline{A}} P_k \phi$.  This gives a function $\Phi_k$ with Fourier support in the region
$\{ \xi: 2^{k-5} \leq |\xi| \leq 2^{k+5} \}$ obeying the bounds
\begin{align}
\| \Phi_k \|_{S} &\lesssim K_k \label{phik-1}\\
\| (\Phi_k - P_k \phi)[t_0] \|_D &\lesssim \eps^{\delta} K_k\label{phik-2}\\
\| \Box'_{P_{\leq k-20} \underline{A}} (\Phi_k - P_k \phi) \|_{N_2} &\lesssim \eps^{\delta} K_k\label{phik-3}
\end{align}
where
$$ K_k := \| P_k \phi[t_0] \|_D + \|  \Box'_{P_{\leq k-20} \underline{A}} P_k \phi \|_{N_2}.$$
We first show the estimate
\begin{equation}\label{eq:kk-bound}
\sum_k K_k^2 \lesssim 1.
\end{equation}
For the $\| P_k \phi[t_0] \|_D$ portion of $K_k$, this follows from the normalization $\| \phi[t_0] \|_D = \| (f,g)\|_D = 1$ and orthogonality.  Thus
it suffices to show that
$$ \sum_k \| \Box'_{P_{\leq k-20} \underline{A}} P_k \phi \|_{N_2}^2 \lesssim 1.$$
Observe that $\Box'_{P_{\leq k-20} \underline{A}} P_k \phi$ has frequency $\sim 2^k$.  Thus we can rewrite the previous as
$$ \sum_k 2^{(n-4)k} \| \Box'_{P_{\leq k-20} \underline{A}} P_k \phi \|_{L^1_t L^2_x}^2 \lesssim 1.$$
Since $P_k \Box'_{\underline{A}} \phi = 0$, it thus suffices to prove the estimates
\begin{equation}\label{eq:comm}
 \sum_k 2^{(n-4)k} \| [\Box'_{P_{\leq k-20} \underline{A}},P_k] \phi \|_{L^1_t L^2_x}^2 \lesssim 1
\end{equation}
and
\begin{equation}\label{eq:comm2}
 \sum_k 2^{(n-4)k} \| P_k (\Box'_{\underline{A}} - \Box'_{P_{\leq k-20} \underline{A}}) \phi \|_{L^1_t L^2_x}^2 \lesssim 1.
\end{equation}
We first prove \eqref{eq:comm}.  The top order term $\Box$ of $\Box'_{P_{\leq k-20} \underline{A}}$ commutes with $P_k$, as does $\nabla_x$, so it suffices to show that
$$
 \sum_k 2^{(n-4)k} \| [P_{\leq k-20} \underline{A}, P_k] \cdot \nabla_x \phi \|_{L^1_t L^2_x}^2 \lesssim 1.$$
We may freely insert a $P_{k-5 < \cdot < k+5}$ in front of $\phi$.
By Lemma \ref{cook-lemma} we may estimate the left-hand side by
$$
\lesssim \sum_k 2^{(n-6)k} \| \nabla_x P_{\leq k-20} \underline{A} \|_{L^2_t L^{2n/3}_x}^2  \| \nabla_x P_{k-5 < \cdot < k+5} \phi \|_{L^2_t L^{2n/(n-3)}_x}^2.$$
However, by \eqref{eq:a-s} we have
\begin{align*}
\| \nabla_x P_{\leq k-20} \underline{A} \|_{L^2_t L^{2n/3}_x}
&\lesssim 
\| \nabla_x \underline{A} \|_{L^2_t \dot B^{[p_*,2n/3]}_2} \\
&\lesssim \eps
\end{align*}
so we can bound the previous expression by
$$ \lesssim \eps^2 
\sum_k 2^{(n-6)k} \| \nabla_x P_{k-5 < \cdot < k+5} \phi \|_{L^2_t L^{2n/(n-3)}_x}^2$$
which in turn is
$$ \lesssim \eps^2 \| \nabla_x \phi \|_{L^2_t \dot B^{[2n/(n-3),2n/3]}_2}^2
\lesssim \eps^2 \| \phi \|_S^2 \lesssim 1$$
as desired, by \eqref{eq:phis-big}.

Now we prove \eqref{eq:comm2}.  Since
\begin{equation}\label{eq:box-box}
 \Box'_{\underline{A}} - \Box'_{P_{\leq k-20} \underline{A}} 
= 2i (P_{>k-20} \underline{A}) \cdot \nabla_x,
\end{equation}
it suffices to show that
$$\sum_k 2^{(n-4)k} \| P_k ((P_{>k-20} \underline{A}) \cdot \nabla_x \phi) \|_{L^1_t L^2_x}^2 \lesssim 1.$$
Discarding the $P_k$, we use H\"older and bound the left-hand side by
$$\lesssim \sum_k 2^{(n-4)k} \| P_{>k-20} \underline{A} \|_{L^2_t L^{2n/(n-3)}_x}^2 \| \nabla_x \phi \|_{L^2_t L^{2n/3}_x}^2.$$
By \eqref{eq:phis-big} we have $\| \nabla_x \phi \|_{L^2_t L^{2n/3}_x} \lesssim \eps^{-\delta/2}$, so we can bound the previous by
$$\lesssim \eps^{-\delta} \sum_k 2^{(n-4)k} \| P_{>k-20} \underline{A} \|_{L^2_t L^{2n/(n-3)}_x}^2,$$
which by splitting $P_{>k-20}$ and rearranging (exploiting the positivity of $n-4$) can be bounded by
\begin{align*} 
&\lesssim \eps^{-\delta} \sum_k 2^{(n-4)k} \| P_k \underline{A} \|_{L^2_t L^{2n/(n-3)}_x}^2 \\
&\lesssim \eps^{-\delta} \sum_k 2^{(n-6)k} \| P_k \nabla_x \underline{A} \|_{L^2_t L^{2n/(n-3)}_x}^2 \\
&\lesssim \eps^{-\delta} \| \nabla_x \underline{A} \|_{L^2_t \dot B^{[2n/(n-3),2n/3]}_2}^2 \\
&\lesssim \eps^{-\delta} \eps^2,
\end{align*}
which is acceptable.  This concludes the proof of \eqref{eq:comm2} and hence \eqref{eq:kk-bound}.

Now we return to the proof of \eqref{phis:targ}.  Since $S$ is built out of $l^2 $ Besov spaces, with all $L^p$ norms greater than or equal to 2, we have
$$ \| \phi \|_S \lesssim (\sum_k \| P_k \phi \|_S^2)^{1/2},$$
and so it suffices to show that
$$ \sum_k \| P_k \phi \|_S^2 \lesssim 1.$$
From \eqref{phik-1}, \eqref{eq:kk-bound} we have 
$$ \sum_k \| \Phi_k \|_S^2 \lesssim 1,$$
so it suffices to show that
$$ \sum_k \| P_k \phi - \Phi_k \|_S^2 \lesssim 1.$$
By \eqref{strichartz-covariant-bootstrap} we have
$$ \| P_k \phi - \Phi_k \|_S \lesssim \eps^{-\delta/2}  (\| (P_k \phi - \Phi_k)[t_0] \|_D + \| \Box'_{\underline{A}} (P_k \phi - \Phi_k) \|_{N_2}),$$
and hence by \eqref{phik-2}, \eqref{phik-3}
$$ \| P_k \phi - \Phi_k \|_S \lesssim \eps^{-\delta/2} (\eps^{\delta} K_k + \| (\Box'_{\underline{A}} - \Box'_{P_{\leq k-20} \underline{A}}) (P_k \phi - \Phi_k) \|_{N_2}).$$
The $\eps^{\delta} K_k$ term is acceptable by \eqref{eq:kk-bound}.  Thus it remains to show that
$$ \sum_k \| (\Box'_{\underline{A}} - \Box'_{P_{\leq k-20} \underline{A}}) (P_k \phi - \Phi_k) \|_{N_2}^2 \lesssim 
\eps^\delta.$$
By \eqref{eq:box-box} it suffices to show that
$$ \sum_k \| (P_{>k-20} \underline{A}) \cdot \nabla_x (P_k \phi - \Phi_k) \|_{N_2}^2 \lesssim \eps^\delta.$$
We rewrite the left-hand side as
$$ \sum_k \| \sum_{m > -20} (P_{k+m} \underline{A}) \cdot \nabla_x (P_k \phi - \Phi_k) \|_{N_2}^2$$
and note that the terms in the summation become almost orthogonal once $m > 20$, and so we can estimate this expression by
$$ \lesssim \sum_k \big (\|\sum_{m > -20}  
\|(P_{k+m} \underline{A}) \cdot \nabla_x (P_k \phi - \Phi_k)\|^2_{\dot H^{\frac n2-2}} \big )^{\frac 12}\|_{L^1_t}^2.$$
Since the expression inside the norm has frequency at most $O(2^{k+m})$, we can estimate the left-hand side by
$$ \lesssim \sum_k 2^{(n-4)(k+m)} \| \big (\sum_{m>-20}
\|(P_{k+m} \underline{A}) \cdot \nabla_x (P_k \phi - \Phi_k) \|_{L^2_x}^2\big )^{\frac 12} \|_{L^1_t}^2,$$
which by H\"older is bounded by
$$ \lesssim \sum_k \|\big (\sum_{m > -20} 2^{(n-4)(k+m)} \| P_{k+m} \underline{A} \|_{L^{2n/(n-3)}_x}^2 \big )^{\frac 12}
\| \nabla_x (P_k \phi - \Phi_k) \|_{L^{2n/3}_x} \|_{L^1_t}^2.$$
This in turn is bounded by
$$ \lesssim \sum_k \sum_{m> -20} 2^{(n-6)(k+m)} \| P_{k+m} \nabla_x \underline{A} \|_{L^2_t L^{2n/(n-3)}_x}^2 \| \nabla_x (P_k \phi - \Phi_k) \|_{L^2_t L^{2n/3}_x}^2.$$
From \eqref{eq:a-s} we see that 
$$ \sum_{m > -20} 2^{(n-6)(k+m)} \| P_{k+m} \nabla_x \underline{A} \|_{L^2_t L^{2n/(n-3)}_x}^2 \lesssim \eps^2,$$
so we can bound the previous by
$$ \lesssim \eps^2 \sum_k \| \nabla_x (P_k \phi - \Phi_k) \|_{L^2_t L^{2n/3}_x}^2.$$
By the triangle inequality, this is bounded by
$$\lesssim \eps^2 \sum_k \| \nabla_x P_k \phi \|_{L^2_t L^{2n/3}_x}^2
+ \eps^2 \sum_k \| \nabla_x \Phi_k \|_{L^2_t L^{2n/3}_x}^2,
$$
which is bounded in turn by
$$ \lesssim \eps^2 \| \phi \|_S^2 + \eps^2 \sum_k \| \Phi_k \|_S^2,$$
and the claim then follows from \eqref{eq:phis-big}, \eqref{phik-1}, \eqref{eq:kk-bound}.
This concludes the proof of Proposition \ref{covariant-Strichartz}, provided that we can construct the parametrix in Proposition \ref{localized-Strichartz}.

To complete the proof of Theorem \ref{main} it thus remains to prove Proposition \ref{localized-Strichartz}.  To begin with, we observe from Duhamel's principle \eqref{duhamel-2} and Minkowski's inequality that it suffices to do this when $F=0$.  Secondly, by scale invariance we may take $k=0$, while from time translation invariance we take $t_0 = 0$.  Thirdly, we may normalize $K = 1$.  Finally, we can take advantage of the frequency localization near the frequency 1 to replace all Besov norms with their Lebesgue counterparts.  We are thus reduced to proving the following proposition.

\begin{proposition}\label{Strichartz-0}  Let $I$ be a compact time interval containing $0$, and let $\underline{A}$ be a smooth connection on $I \times \R^n$ which obeys the smallness condition
\begin{equation}\label{a-small}
\| \underline{A} \|_D + \| \Box \underline{A} \|_{N_1} \lesssim \eps
\end{equation}
on $I \times \R^n$.  We assume that $\underline{A}$ is in the Coulomb gauge \eqref{CG-underline}.
Suppose also that $\underline{A}$ has Fourier support on the region $\{ \xi: |\xi| \leq 2^{-10} \}$.  Let $(f,g)$ be a pair of Schwartz functions on $\R^n$ with Fourier support in the region $\{ \xi: 2^{-3} \leq |\xi| \leq 2^{3} \}$, and with norm
\begin{equation}\label{eq:fg-norm}
 \| f \|_2 + \|g\|_2 = 1.
\end{equation}
Then there exists a function $\phi$ on $I \times \R^n$ 
with Fourier support in the region $\{ \xi: 2^{-10} \leq |\xi| \leq 2^{10} \}$ obeying the estimates 
\begin{equation}\label{strichartz-covariant-localized}
\| \nabla_{x,t} \phi \|_{L^\infty_t L^2_x} + \| \nabla_{x,t} \phi \|_{L^2_t L^{p_*}_x} \lesssim 1,
\end{equation}
\begin{equation}\label{eq:data-approx}
\| \phi(0) - f \|_{L^2_x} + \| \phi_t(0) - g \|_{L^2_x} \lesssim \eps^{\delta},
\end{equation}
and
\begin{equation}\label{eq:forcing-approx}
\| \Box'_{\underline{A}} \phi \|_{L^1_t L^2_x} \lesssim \eps^{\delta}
\end{equation}
on $I \times \R^n$.
\end{proposition}

We shall prove this Proposition in the remainder of this paper.  For now, we content ourselves with making more two reductions for the above Proposition.  First, we remark that the condition that $\phi$ have Fourier support in the region $\{ \xi: 2^{-10} \leq |\xi| \leq 2^{10} \}$ can be dropped.  For, if we can find another function $\tilde \phi$ on $I \times \R^n$ which already obeys \eqref{strichartz-covariant-localized}, \eqref{eq:data-approx}, \eqref{eq:forcing-approx} but without the frequency support assumption, then we claim that the function $\phi := P_{-5 \leq \cdot \leq 5 } \tilde \phi$ will also obey these three estimates while also having the Fourier support property.  The verification of this for \eqref{strichartz-covariant-localized} and \eqref{eq:data-approx} is immediate, since we can just apply $P_{-5 \leq \cdot \leq 5}$ to the left-hand side (and use the Fourier support of $(f,g)$).  Now we verify \eqref{eq:forcing-approx}.  Since we already assume
$$ 
\| \Box'_{\underline{A}} \tilde \phi \|_{L^1_t L^2_x} \lesssim \eps^{\delta}$$
we thus have
$$ 
\| P_{-5 \leq \cdot \leq 5} \Box'_{\underline{A}} \tilde \phi \|_{L^1_t L^2_x} \lesssim \eps^{\delta}$$
and thus we only need to prove the commutator estimate
$$ 
\| [P_{-5 \leq \cdot \leq 5}, \Box'_{\underline{A}}] \tilde \phi \|_{L^1_t L^2_x} \lesssim \eps^{\delta}.$$
Since the top order term $\Box$ commutes with $P_{-5 \leq \cdot \leq 5}$, as does $\nabla_x$, it suffices to show
$$ 
\| [P_{-5 \leq \cdot \leq 5}, \underline{A}]  \cdot \nabla_x \tilde \phi \|_{L^1_t L^2_x} \lesssim \eps^{\delta}.$$
By Lemma \ref{cook} we can estimate the left-hand side by
$$ \| \nabla_x \underline{A} \|_{L^2_t L^q_x}  \| \nabla_x \tilde \phi \|_{L^2_t L^{p_*}_x}$$
where $1/q + 1/p_* = 1/2$.  However, by definition of $p_*$, we have $q > 2n/3$, and so by Bernstein (using the frequency localization of $\underline{A}$) and \eqref{eq:a-s} we have
$$ \| \nabla_x \underline{A} \|_{L^2_t L^q_x} 
\lesssim \| \underline{A} \|_S \lesssim \eps.$$

The claim then follows from \eqref{strichartz-covariant-localized} (for $\tilde \phi$).  

Our next remark concerns the spacetime Fourier support of $\underline{A}$.  We may extend $\underline{A}$ from $I \times \R^n$ to all of $\R \times \R^n$ by evolving by the free wave equation both forward and backward in time.  By hypothesis, $\underline{A}$ then obeys the Cauchy problem
\begin{equation}\label{eq:cauchy}
\begin{split}
\Box \underline{A} &= F \\
\underline{A}(0) &= f \\
\partial_t \underline{A}(0) &= g
\end{split}
\end{equation}
where $F$, $f$, $g$ are divergence-free and have spatial Fourier support on the region $\{ |\xi| \leq 2^{-10} \}$, and we have the estimates
$$ \| f \|_{\dot H^{n/2-1}_x} + \| g \|_{\dot H^{n/2-2}_x} + \| F \|_{L^1_t \dot B^{[2,n/2]}_1} \lesssim \eps.$$
Now consider the spacetime $\tilde F(\tau,\xi)$ of the Fourier transform of $F$.  We claim that we may restrict this spacetime Fourier transform of $F$ to the region\footnote{Of course, this will destroy any compact support properties that $F$ has in time, but this will not affect the rest of the argument.}
\begin{equation}\label{eq:sector}
 \{ (\tau,\xi): |\xi| \leq 2^{-10}; |\tau| \leq 4|\xi| \}.
\end{equation}
Indeed, if $\tilde F$ does not lie in this region, we can smoothly decompose $F = F_0 + \sum_{k > 0} F_k$, where $\tilde F_0$ has the correct Fourier support, and $\tilde F_k$ has Fourier support on the region
$$ \{ (\tau,\xi): |\xi| \leq 2^{-10}; 2^{k-1} |\xi| \leq |\tau| \leq 2^{k+1} |\xi| \}.$$
Some Littlewood-Paley theory reveals that all functions remain divergence-free, and we have the bounds
$$ \| F_k \|_{L^1_t \dot B^{[2,n/2]}_1} \lesssim \| F \|_{L^1_t \dot B^{[2,n/2]}_1} \lesssim \eps$$
for both $k=0$ and $k>0$.  We can then split $\underline{A} = \underline{A_0} + \sum_k \frac{F_k}{\Box}$, where $\underline{A_0}$ solves the Cauchy problem
\begin{align*}
\Box \underline{A_0} &= F_0 \\
\underline{A_0}(0) &= f \\
\partial_t \underline{A_0}(0) &= g
\end{align*}
and $\frac{1}{\Box}$ is the spacetime Fourier multiplier that inverts $\Box$ (this is well-defined on $F_k$ since the Fourier support of $F_k$ avoids the light cone).  Note that $\underline{A_0}$ remains divergence-free.

We observe that on the Fourier support of $F_k$, the spacetime Fourier multiplier $\frac{2^{2k} \Delta}{\Box}$ has a bounded smooth symbol, so we have
$$ \| \frac{F_k}{\Box} \|_{L^1_t \dot B^{[2,\infty]}_1} \lesssim 2^{-2k} \| F_k \|_{L^1_t \dot B^{[2,n/2]}_1} \lesssim 2^{-2k} \eps.$$
In particular we have
$$ \| \sum_{k > 0} \frac{F_k}{\Box} \|_{L^1_t L^\infty_x} \lesssim 
\| \sum_{k > 0} \frac{F_k}{\Box} \|_{L^1_t \dot B^{[2,\infty]}_1} \lesssim \eps$$
(note how the $l^1$ nature of the Besov norm $\dot B^{[2,\infty]}_1$ is needed here).  Because of this $L^1_t L^\infty_x$ nature of $\underline{A} - \underline{A_0} = \sum_{k > 0} \frac{F_k}{\Box}$, any parametrix $\phi$ which obeys the properties \eqref{strichartz-covariant-localized}, \eqref{eq:data-approx}, \eqref{eq:forcing-approx} for $\Box_{\underline{A_0}}$ will also obey the same properties for $\Box_{\underline{A}}$.  Indeed, we have
$$ \| (\Box_{\underline{A}} - \Box_{\underline{A_0}}) \phi \|_{L^1_t L^2_x}
\lesssim \| \sum_{k > 0} \frac{F_k}{\Box} \cdot \nabla_x \phi \|_{L^1_t L^2_x}
\lesssim \| \nabla_x \phi \|_{L^\infty_t L^2_x} \lesssim \eps^\delta.$$

Thus we may assume without loss of generality that $F$ has spacetime Fourier support in the region \eqref{eq:sector}.  In particular, this gives us good control on the time regularity of $F$ (and hence of $A$).  In particular, we see that time derivatives behave the same way as spatial derivatives in our estimates on $A$ and $F$, for instance we have
\begin{equation}\label{eq:f-bound} \| \nabla_{x,t} F \|_{L^1_t \dot B^{[2,n/3]}_1} \lesssim \eps.
\end{equation}
In particular, by integrating $F$ in time we obtain
\begin{equation}\label{eq:f-sup}
 \| F \|_{L^\infty_t \dot B^{[2,n/3]}_1} \lesssim \eps.
\end{equation}

\section{Distorted plane waves}\label{sec:distorted}

We now begin the proof of Proposition \ref{Strichartz-0}.  Fix the (small, frequency-localized) divergence-free connection $\underline{A}$; we will assume that the forcing term $F := \Box \underline{A}$ in this connection is defined on $\R \times \R^n$ and has spacetime Fourier support in \eqref{eq:sector}.  We are given frequency-localized initial data $f,g$, and wish to construct an approximate solution $\phi$ to the equation $\Box'_{\underline{A}} \phi = 0$ with initial data $\phi[0] = (f,g)$, and which obeys the Strichartz estimates.  As remarked in the previous section, we do not need to enforce any frequency localization properties on $\phi$.

We shall tackle this problem in three stages.  Firstly, by constructing a family of distorted plane waves, we shall build a large class of functions $\phi$ which obey Strichartz and energy estimates.  Secondly, we use our energy estimates to show that for our given initial data $(f,g)$ there exists a member $\phi$ of this class which is close to $(f,g)$ at time 0.  Finally, we show that these class of functions are good approximate solutions to the covariant wave equation $\Box'_{\underline{A}} \phi = 0$.  
Although, the construction of the distorted plane waves will not rely on the fact that
the connection ${\underline{A}} $ is in the Coulomb gauge \eqref{CG-underline},
the latter condition becomes important in lowering the dimension $n$ in  Proposition \ref{Strichartz-0} to $n\ge \n$.

We begin with the construction of the distorted plane waves.  We first give an informal discussion.  For the free wave equation, we recall that the plane waves
$$ e^{2\pi ix \cdot \xi} e^{\pm 2\pi i t |\xi|}$$
are exact solutions to the free wave equation $\Box \phi = 0$ for any $\xi \in \R^n$ and any sign $\pm$.  In particular, we can construct frequency-localized free waves by the representation
$$ \phi(t,x) := \int e^{2\pi i x \cdot \xi} e^{\pm 2\pi i t|\xi|} h_\pm(\xi) a(\xi)\ d\xi$$
where $h_\pm$ is an arbitrary (smooth) function and $a(\xi)$ is a fixed cutoff to the region $2^{-10} \leq |\xi| \leq 2^{10}$ which equals 1 on $2^{-5} \leq |\xi| \leq 2^5$.  The functions $h_\pm$ are essentially the Fourier transforms of the initial data $\phi[0] = (\phi(0), \phi_t(0))$, and we have free Strichartz estimates such as
$$
\| \nabla_{x,t} \phi \|_{L^\infty_t L^2_x} + \| \nabla_{x,t} \phi \|_{L^2_t L^{p_*}_x} \lesssim \| h_\pm \|_2,$$
from \eqref{strichartz} and Plancherel.

Motivated by this, we now construct a similar parametrix for the modified D'Alambertian $\Box'_{\underline{A}} = \Box + 2 i \underline{A} \cdot \nabla_x$.  Just as the ODE operator $\partial_x + iA$ can be conjugated to $\partial_x$ by means of the integrating factor $\exp(i \partial_x^{-1} A)$, we would expect that $\Box'_{\underline{A}}$ should be somehow approximately conjugate to $\Box$ using some phase correction roughly of the form $\exp(i \nabla^{-1} \underline{A})$.  Accordingly, we shall consider waves of the form
$$ \phi(t,x) := U_\pm(t) h_\pm$$
where $\pm$ is a sign, $h_\pm(\xi)$ is an arbitrary smooth function, $U_\pm(t)$ is the operator
\begin{equation}\label{eq:ut-def}
 U_\pm(t) h_\pm(x) := \int e^{2\pi i\Psi_\pm(t,x,\xi)} e^{2\pi i x \cdot \xi} e^{\pm 2\pi i t|\xi|} h_\pm(\xi) a(\xi)\ d\xi 
\end{equation}
and $\Psi_\pm(t,x,\xi)$ is a smooth, real-valued phase correction, to be chosen later, which should heuristically have the scaling and frequency of $\nabla^{-1} \underline{A}$.

To motivate how we shall choose $\Psi_\pm$, let us compute
$$ \Box'_{\underline{A}} \phi(t,x) = (-\partial_t^2 + \Delta + 2 i \underline{A} \cdot \nabla_x) (U_\pm(t) h_\pm).$$
A brief calculation shows that this is equal to
\begin{equation}\label{eq:ut-box}
\Box'_{\underline{A}} \phi(t,x) = 2\pi \int \Omega_\pm(t,x,\xi)
e^{2\pi i\Psi_\pm(t,x,\xi)} e^{2\pi i x \cdot \xi} e^{\pm 2\pi i t|\xi|} h_\pm(\xi) a(\xi)\ d\xi 
\end{equation}
where $\Omega_\pm$ is the function
\begin{equation}\label{eq:omega-def}
 \Omega_\pm := -4\pi |\xi| L^\mp_\omega \Psi_\pm - 2 \underline{A} \cdot \xi
+ i \Box \Psi_\pm + 2\pi (|\partial_t \Psi_\pm|^2 - |\nabla_x \Psi_\pm|^2) - 2 \underline{A} \cdot \nabla_x \Psi_\pm
\end{equation}
where $\omega := \xi/|\xi|$ is the direction of $\xi$, and $L^\mp_\omega$ is the null vector field
$$L_\omega^\mp := \omega \cdot \nabla_x \mp \partial_t.$$

Clearly, in order to make $\Box'_{\underline{A}} \phi$ small, it will be helpful to make $\Omega_\pm$ small as well.  The last three terms in \eqref{eq:omega-def} are quadratic in $\underline{A}$ and $\Psi$ (and morally have the scaling of $A^2$) and will be easy to manage.  The term $\Box \Psi_\pm$ will turn out to be very small because\footnote{Actually, this additional structure is unnecessary in high dimensions; the mere fact that $\Psi_\pm$ will look like $\nabla^{-1} \underline{A}$ will make $\Box \Psi_\pm$ have the scaling of $\nabla A$, which is comparable in strength to $A^2$ and will be manageable since the derivative is falling on a low frequency term.  By a similar token, the null structure in $|\partial_t \Psi_\pm|^2 - |\nabla_x \Psi_\pm|^2$ will be ignored in this high-dimensional setting, though it may well play a role in lower dimensions.}
 $\Box \underline{A}$ is small and we will construct $\Psi_\pm$ to look roughly like $\nabla^{-1} \underline{A}$.  

The main terms are the first two.  Neither of them are individually manageable (their contribution has the scaling of $A \nabla \phi$, which is not controllable by Strichartz estimates).  However, we will choose $\Psi_\pm$ so that these terms mostly cancel, i.e. we will choose $\Psi_\pm$ so that\footnote{As we shall see in the ensuing discussion, the magnitude $r = |\xi|$ of the frequency plays almost no role.  Indeed it is possible to take a Fourier transform in the radial variable, and replace the plane wave approximation with a moving plane approximation, based on a (phase-)distorted Radon transform instead of a distorted Fourier transform.  We however will retain the Fourier transform-based approach as it is a more familiar approach for building parametrices.}
$$ 2\pi L_\omega^\mp \Psi_\pm \approx -\underline{A} \cdot \omega.$$
If $L_\omega^\mp$ were elliptic, it would thus make sense to choose\footnote{This choice of $\Psi_\pm$ has a geometric interpretation; it corresponds to the gauge change which would make $A$ vanish in the direction $L_\omega^\mp$; in other words, it is a Cronstrom-type gauge in the specified null direction.  However, because this gauge depends on the choice of null vector field $L_\omega^\mp$, it cannot be represented by a single, global change of gauge for the connection $\underline{A}$.  Thus one can view the phase correction $e^{2\pi i \Psi_\pm}$ in the parametrix as a sort of ``microlocal Cronstrom gauge''.}
$$ \Psi_\pm := -\frac{1}{2\pi} (L_\omega^\mp)^{-1} \underline{A} \cdot \omega.$$
Unfortunately, the differential operator $L_\omega^\mp$ is not elliptic, indeed its symbol vanishes on the null plane $\Pi^\mp_\omega := \{ (\tau,\xi): \xi \cdot \omega = \pm \tau\}$.  Fortunately, we are assuming that $\underline{A}$ is almost a solution to the free wave equation $\Box \underline{A} = 0$, which means, morally speaking, that $\underline{A}$ has Fourier support on the light cone $\{ (\tau,\xi): |\xi| = |\tau|\}$.  Because the null plane $\Pi^\mp_\omega$ is tangent to this light cone, the symbol $L_\omega^\mp$ behaves as if it is elliptic on free waves, and this insight shall be key in making this parametrix function effectively.  (The forcing term $\Box \underline{A}$ will cause some technical difficulties, but they will be manageable because we insisted that we have an $l^1$-Besov control on $\Box \underline{A}$, and not just the more familiar $L^1_t \dot H^{n/2-2}_x$ control.)  We remark that for technical reasons we shall have to truncate $\underline{A}$ in frequency space near the null plane $\Pi^\mp_\omega$ to avoid singularities, and in particular to keep $\Psi$ smooth in the $\omega$ variable.  This will unfortunately introduce a number of additional error terms
into our analysis; in very high dimensions (e.g. $n \geq 10$) these error terms can be relatively easily handled
by choosing the truncation parameter $\sigma$ appropriately, but in order to handle the medium 
dimensions $6 \leq n \leq 9$ we will require a somewhat delicate analysis exploiting the Coulomb gauge property
\eqref{CG-underline}.

\section{Construction of $\Psi_\pm$}\label{sec:construct}

We now construct $\Psi_\pm$ more rigorously.  Let $\omega \in S^{n-1}$ be any unit vector in $\R^n$.  Our starting point is the null frame decomposition
\begin{equation}\label{eq:box-null}
 \Box =  L_\omega^+ L_\omega^- + \Delta_{\omega^\perp}
\end{equation}
of the free D'Alambertian $\Box$, where
$$ \Delta_{\omega^\perp} := \Delta - (\omega \cdot \nabla_x)^2$$
is the Laplacian for the hyperplane in $\R^n$ orthogonal to $\omega$.  In particular, if $\underline{A}$ solves the free wave equation $\Box \underline{A} = 0$, then we have
$$ \underline{A} = -L_\omega^\mp (L_\omega^\pm \Delta_{\omega^\perp}^{-1} \underline{A})$$
(noting that $L_\omega^+$ and $L_\omega^-$ are constant-coefficient differential operators for fixed $\omega$, and thus commute).
Motivated by this, it seems reasonable to choose $\Psi_\pm$ via the formula
$$ \Psi_\pm ``:='' \frac{1}{2\pi} L_\omega^\pm \Delta_{\omega^\perp}^{-1} \underline{A} \cdot \omega.$$
Indeed, if we adopted this definition of $\Psi_\pm$, then we see from \eqref{eq:box-null} that we would have
$$ 2\pi L_\omega^\mp \Psi_\pm +\underline{A} \cdot \omega = 
\Box \Delta_{\omega^\perp}^{-1} \underline{A} \cdot \omega =  \Delta_{\omega^\perp}^{-1} F\cdot \omega $$
where $F := \Box \underline{A}$ is as before.  

This choice of definition for $\Psi_\pm$
has several good properties; not only is $2\pi L_\omega^\mp \Psi_\pm + \underline{A} \cdot \omega$ small, but $\Psi_\pm$ also obeys a number of good bounds, for instance its $L^1_t L^\infty_x$ norm can be proven to be small.  However, the singularity of $\Delta_{\omega^\perp}$ along spatial frequencies parallel to $\omega$ makes this choice of $\Psi_\pm$ very rough with respect to the $\omega$ (and thus $\xi$) variables.  To fix this problem we will smooth $\Psi_\pm$ out near this frequency singularity.

More precisely, for any direction $\omega \in S^{n-1}$ and any angle $0 < \theta \lesssim 1$, we define the sector projection $\Pi_{\omega,>\theta}$ in frequency space by the formula
$$ \widehat{\Pi_{\omega,>\theta} f}(\xi) := (1 - \eta(\frac{\angle(\xi,\omega)}{\theta}))  (1 - \eta(\frac{\angle(-\xi,\omega)}{\theta})) \hat f(\xi) $$
where $\eta(\xi)$ is a bump function on $\R^n$ which equals 1 when $|\xi| < 1/2$ and vanishes for $|\xi| > 1$, and $\angle(\xi,\omega)$ is the angle between $\xi$ and $\omega$.  Thus $\Pi_{\omega, >\theta}$ restricts $f$ smoothly (except at the frequency origin) to the sector of frequencies $\xi$ whose angle with both $\omega$ and $-\omega$
is $\gtrsim \theta$.  We remark that this operator $\Pi_{\omega, > \theta}$ is a Fourier multiplier and in particular commutes with Littlewood-Paley projections and constant-coefficient differential operators. It is also important to note that  $\Pi_{\omega, > \theta}$
preserves the space of real-valued functions.
We also define the complementary operator $\Pi_{\omega, \leq \theta}$ by
$$ \Pi_{\omega, \leq \theta} := 1 - \Pi_{\omega, > \theta}.$$

Let $\sigma > 0$ be a small exponent (depending only on $n$) to be chosen later.  
We will define $\Psi_\pm$ by the formula
\begin{equation}\label{eq:psi-pm}
 \Psi_\pm := \frac{1}{2\pi} L_\omega^\pm \Delta_{\omega^\perp}^{-1} \sum_{k < -5} \Pi_{\omega, > 2^{\sigma k}} P_k \underline{A} \cdot \omega,
\end{equation}
where $\omega := \xi/|\xi|$.  (We do not need to define $\Psi_\pm$ at the frequency origin, since our parametrix \eqref{eq:ut-def} is restricted to the region $|\xi| \sim 1$).  Note that if it were not for the projections $\Pi_{\omega, > 2^{\sigma k}}$, then this definition would match the previous proposal for $\Psi_\pm$ (note that the frequency restrictions on $\underline{A}$ ensure that $\underline{A} = \sum_{k < -5} P_k \underline{A}$).
We make the fundamental observation that $\Psi_\pm$ is \emph{real-valued}; this is because $\underline{A}$ is real valued,
and all the Fourier multipliers here have symbols that are real and even.

Now apply $2\pi L_\omega^\mp$ to both sides of \eqref{eq:psi-pm}.  By \eqref{eq:box-null}, we have
$$
2\pi L_\omega^\mp \Psi_\pm := -\sum_{k < -5} \Pi_{\omega, > 2^{\sigma k}} P_k \underline{A} \cdot \omega + \sum_{k < -5} \Delta_{\omega^\perp}^{-1} \Pi_{\omega, > 2^{\sigma k}} P_k F \cdot \omega.$$
In particular, we have
\begin{equation}\label{eq:psi-small}
 2\pi L_\omega^\mp \Psi_\pm +\underline{A} \cdot \omega = 
\sum_{k < -5} \Pi_{\omega, \leq 2^{\sigma k}} P_k \underline{A} \cdot \omega +
\sum_{k < -5} \Delta_{\omega^\perp}^{-1} \Pi_{\omega, > 2^{\sigma k}} P_k F \cdot \omega.
\end{equation}
The first term on the right-hand side of \eqref{eq:psi-small} will generate error terms in our parametrix which resemble the bad expression $\underline{A} \cdot \nabla_x \phi$, but fortunately the additional projection $\Pi_{\omega, \leq 2^{\sigma k}}$ will render these types of expressions manageable, provided that $\sigma$ is sufficiently large.
The second term, meanwhile, will benefit from the Fourier multiplier 
$\Delta_{\omega^\perp}^{-1} \Pi_{\omega, > 2^{\sigma k}}$, which is bounded by $O(2^{-2\sigma k})$.  In fact
we will be able to exploit the divergence-free nature of $F$ (and the dot product with $\omega$) to effectively
improve this estimate further, basically to $O(2^{-3\sigma k})$).

Observe that the function $\Psi_\pm(t,x,\xi)$ depends only on the angular component $\omega := \xi/|\xi|$ of the frequency $\xi$, and not on the radial component $|\xi|$.  In particular, $\Psi_\pm$ is infinitely smooth (indeed, it is constant) in the radial direction $\xi \cdot \nabla_\xi$.  The regularity in the $\omega$ direction is not as good, however we will be able to obtain satisfactory control this direction if $\sigma$ is sufficiently small (since we will then be truncating to be well away from the singularities of $\Delta_{\omega^\perp}^{-1}$).  Thus in order to keep both error terms in \eqref{eq:psi-small} under control, $\sigma$ can neither be too large or too small.
Specifically, we will choose $\sigma$ so that
\begin{equation}\label{sigma-bounds}
0<\frac{n+1}{(n-1)(n-3)} < \sigma <  \frac{1}{2}.
\end{equation}
One can easily verify that such a $\sigma$ can be chosen whenever $n\geq \n$.

\section{Interlude: decomposable functions}\label{sec:interlude}

We now pause to discuss an abstract problem concerning how to deal with certain factors embedded in the kernel of linear operators.  

Compare the expression $U^\pm(t) h$ defined in \eqref{eq:ut-def}, with the expression $\Box'_{\underline{A}} U^\pm(t) h$ in \eqref{eq:ut-box}.  The two expressions are almost the same except for an additional factor $\Omega_\pm(t,x,\xi)$ in the second expression.  Heuristically, this means that if one has good estimates on \eqref{eq:ut-def} and good estimates on $\Omega_\pm$, then this should automatically imply good estimates for \eqref{eq:ut-box}, perhaps by some application of H\"older's inequality (or some other bilinear estimate).

Unfortunately one cannot do this directly because the multiplier $\Omega_\pm$ depends on $\xi$ and so cannot be easily pulled outside of the integral.  Nevertheless, as we shall see, the dependence on $\xi$ is sufficiently mild that we will be able to ``decouple'' this multiplier $\Omega_\pm$ from the operator $U^\pm(t)$.

To make this more precise we need some notation.  Let $\Sigma$ be the annulus $\Sigma := \{ \xi: 2^{-10} \leq |\xi| \leq 2^{10} \}$.

\begin{definition}\label{def:multiplier}  Let $1 \leq q,r \leq \infty$, and let $F(t,x,\xi)$ be a function on $\R \times \R^n \times \Sigma$.  We say that $F$ is \emph{decomposable in $L^q_t L^r_x$} if there exists a constant $M > 0$ with the following property: whenever $1 \leq q_1, r_1, q_2, r_2 \leq \infty$ are such that $\frac{1}{q_2} = \frac{1}{q} + \frac{1}{q_1}$ and $\frac{1}{r_2} = \frac{1}{r} + \frac{1}{r_1}$, and whenever $K(t,x,\xi)$ is any kernel for which one has the estimate
$$ \| \int_\Sigma K(t,x,\xi) h(\xi)\ d\xi \|_{L^{q_1}_t L^{r_1}_x} \leq B \| h \|_{L^2_\xi}$$
for all $h \in L^2_\xi$ and some constant $B > 0$, then one must also have the estimate
$$ \| \int_\Sigma F(t,x,\xi) K(t,x,\xi) h(\xi)\ d\xi \|_{L^{q_2}_t L^{r_2}_x} \leq M B \| h \|_{L^2_\xi}$$
for all $h \in L^2_\xi$.  We define $\|F\|_{D(L^q_t L^r_x)}$ to be the the infimum of all constants $M$ with the above property.
\end{definition}

For instance, if $F(t,x,\xi) = F(t,x)$ does not depend on $\xi$, and lies in $L^q_t L^r_x$, then it is automatically decomposable in $L^q_t L^r_x$, and one has
$$ \| F \|_{D(L^q_t L^r_x)} \leq \| F \|_{L^q_t L^r_x};$$
this is because one can simply pull $F(t,x)$ outside of the $\xi$ integral and then apply H\"older's inequality.  Conversely, for general $F(t,x,\xi)$ we have
\begin{equation}\label{eq:qr}
\sup_{\xi \in \Sigma} \| F(\xi) \|_{L^q_t L^r_x} \leq \| F \|_{D(L^q_t L^r_x)},
\end{equation}
as can be seen by specializing to the case $K(t,x,\xi) := K(\xi)$ and $q_1=r_1=\infty$, where $K$ ranges
over arbitrary $L^2_\xi$ functions.

At the other extreme, if $F(t,x,\xi) = F(\xi)$ does not depend on $t,x$, and is bounded, then it is automatically decomposable in $L^\infty_t L^\infty_x$, and one has 
$$ \| F \|_{D(L^\infty_t L^\infty_x)} \leq \| F \|_{L^\infty_\xi(\Sigma)},$$
since one can simply absorb $F$ into the $h$ factor.

It is also easy to see that the $\| \|_{D(L^q_t L^r_x)}$ norm is indeed a norm, and is invariant under translations in time and space.  We also have the algebra property
\begin{equation}\label{eq:algebra}
 \| FG \|_{D(L^q_t L^r_x)} \leq \| F \|_{D(L^{q_1}_t L^{r_1}_x)} \| G \|_{D(L^{q_2}_t L^{r_2}_x)}
\end{equation}
whenever $1/q = 1/q_1 + 1/q_2$ and $1/r = 1/r_1 + 1/r_2$; this is a simple consequence of iterating the above definition.  In particular we have
\begin{equation}\label{eq:factor}
 \| F g(\xi) \|_{D(L^q_t L^r_x)} \leq \| F \|_{D(L^q_t L^r_x)} \|g\|_{L^\infty_\xi(\Sigma)}
\end{equation}
whenever $g$ is independent of $t$ and $x$.

If $F(t,x,\xi)$ is decomposable in $L^q_t L^r_x$, then one should think of this heuristically as if the $\xi$ dependence in $F$ could be ignored, and $F$ could be treated as an $L^q_t L^r_x$ function, at least for the purpose of proving estimates of the above form.

We now show that functions which have enough smoothness in the $\xi$ variable at some angular scale $0 < \theta \lesssim 1$
are also decomposable (uniformly in the choice of $\theta$).

\begin{lemma}\label{banach}  Let $1 \leq q,r \leq \infty$ and $0 < \theta \lesssim 1$.  Let $F(t,x,\xi)$ be a function which is homogeneous in $\xi$ of degree 0 (i.e. $F(t,x,\lambda \xi) = F(t,x,\xi)$ for all $\lambda > 0$).  Then  we have 
$$ \| F \|_{D(L^q_t L^r_x)}
\lesssim \sum_{l=0}^{100n} 
(\theta^{1-n} \int_{\Sigma} \| (\theta \nabla_\xi)^l F(t,x,\xi) \|_{L^q_t L^r_x}^2\ d\xi)^{1/2}$$
\end{lemma}

\begin{proof}
Let us first prove this lemma under the additional assumption that the $\xi$-support of $F$ is contained in a tube of the form $\{ \xi \in \Sigma: \angle(\xi, \omega_0) \lesssim \theta \}$ for some 
$\omega_0 \in S^{n-1}$.  Without loss of generality we may assume that $\omega_0$ is equal to $e_1$, the first basis vector.  Thus $F(t,x,\xi)$ is supported on the tube $T$ where $|\xi| \sim 1$ and $\xi_2, \ldots, \xi_n = O(\theta)$.  In particular, we may write $F = F \psi(\xi)$ where $\psi(\xi)$ is a bump function which equals 1 on the above tube $T$ and is adapted to a slight dilate of that tube.

Write $K := \sum_{l=0}^{100n} (\theta^{1-n} \int_{\Sigma} \| (\theta \nabla_\xi)^l F(t,x,\xi) \|_{L^q_t L^r_x}^2\ d\xi)^{1/2}$.  By Sobolev inequality and homogeneity 
of $F$ we see that $\| (\theta \nabla_\xi)^l F(t,x,\xi) \|_{L^q_t L^r_x} \lesssim K$ for all $\xi \in T$ and all $0 \leq l \leq 50n$.  In particular if we write $F$ as a Fourier series on the support of $\psi$,
$$ F(t,x,\xi) = \sum_{(k_1,\ldots,k_n) \in \Z^n} c_{k_1,\ldots,k_n}(t,x) e^{i(k_1 \xi_1 + \frac{1}{\theta} (k_2 \xi_2 + \ldots + k_n \xi_n))/C}$$
for all $\xi$ in the support of $\psi$, where $C$ is a large constant, then we see from the Fourier inversion formula, integration by parts, and homogeneity of $F$,  that $\|c_k\|_{L^q_t L^r_x} \lesssim \langle k \rangle^{-10n} K$.
Since $F = F\psi$, we thus have
$$ F(t,x,\xi) = \sum_{(k_1,\ldots,k_n) \in \Z^n} c_{k_1,\ldots,k_n}(t,x) e^{i(k_1 \xi_1 + \frac{1}{\theta} (k_2 \xi_2 + \ldots + k_n \xi_n))/C}\psi(\xi)$$
for all $\xi \in \Sigma$.  In particular, from the triangle inequality and \eqref{eq:factor}, we see that $\| F\|_{D(L^q_t L^r_x)} \lesssim K$ as desired.

Now we prove the general case, where $F$ is not supported in a sector.  In this case we can decompose $F = \sum_T F_T$, where $T$ ranges over a finitely overlapping collection of tubes of the
form $T := \{ \xi \in \Sigma: \angle(\xi, \omega_T) \lesssim \theta \}$ for some direction $\omega_T$, and $F_T$ is $F$ smoothly cut off to $T$.
Let $K, q_1, q_2, r_1, r_2, B$ be as in Definition \ref{def:multiplier}.  We observe that
\begin{align*}
\| \int_\Sigma F(t,x,\xi) K(t,x,\xi) h(\xi)\ d\xi \|_{L^{q_2}_t L^{r_2}_x} &\leq \sum_T \| \int_\Sigma F_T(t,x,\xi) K(t,x,\xi) h(\xi)\ d\xi \|_{L^{q_2}_t L^{r_2}_x} \\
&\leq \sum_T \| \int_\Sigma F_T(t,x,\xi) K(t,x,\xi) \chi_T h(\xi)\ d\xi \|_{L^{q_2}_t L^{r_2}_x} \\
&\leq \sum_T \| F_T \|_{D(L^q_t L^r_x)} B \| \chi_T h \|_{L^2_\xi}\\
&\leq B (\sum_T \| F_T \|_{D(L^q_t L^r_x)}^2)^{1/2}
(\sum_T \| \chi_T h \|_{L^2_\xi}^2)^{1/2}\\
&\leq B (\sum_T \sum_{l=0}^{100n} 
\theta^{1-n} \int_{\Sigma} \| (\theta \nabla_\xi)^l F_T(t,x,\xi) \|_{L^q_t L^r_x}^2\ d\xi)^{1/2}\\ &\quad \times (\sum_T \| \chi_T h \|_{L^2_\xi}^2)^{1/2}\\
&\lesssim B (\sum_{l=0}^{100n} 
\theta^{1-n} \int_{\Sigma} \| (\theta \nabla_\xi)^l F(t,x,\xi) \|_{L^q_t L^r_x}^2\ d\xi)^{1/2} \| h \|_{L^2_\xi}
\end{align*}
where in the last line we have used the finite overlap of the $T$.  The claim then follows.
\end{proof}

\section{Estimates on the phase correction $\Psi_\pm$}\label{sec:psi}

We now prove some key estimates on $\Psi_\pm$ and the associated function $\Omega_\pm$, showing that they
can be adequately controlled in the $D(L^q_t L^r_x)$ norms defined earlier and thus these factors can be easily 
dealt with when we prove the Strichartz and energy type estimates.  To obtain these estimates we will first need to decompose the projection $\Pi_{\omega, > 2^{\sigma k}}$ into dyadic pieces.  More precisely, we have
$$ \Pi_{\omega, > 2^{\sigma k}} = \sum_{j=1}^\infty \Pi_{\omega, 2^{\sigma k + j}}$$
where
$$ \Pi_{\omega, \theta} := \Pi_{\omega, > \theta/2} - \Pi_{\omega, > \theta}.$$
In particular, we have
\begin{equation}\label{eq:psi-dyadic}
\Psi_\pm = \frac{1}{2\pi} L_\omega^\pm \Delta_{\omega^\perp}^{-1} \sum_{k < -5} \sum_{j \geq 1} \Pi_{\omega, 2^{\sigma k + j}} P_k \underline{A} \cdot \omega.
\end{equation}
Since $L_\omega^\pm$ can be rewritten as $(\omega,\pm) \cdot \nabla_{x,t}$, we can rewrite the above as
\begin{equation}\label{eq:psi-dyadic-alt}
\Psi_\pm = c_\pm(\omega) \nabla_{x,t} \Delta_{\omega^\perp}^{-1} \sum_{k < -5} \sum_{j \geq 1} \Pi_{\omega, 2^{\sigma k + j}} P_k \underline{A}\cdot \omega
\end{equation}
where $c_\pm(\omega)$ is a smooth, tensor-valued function of $\omega$ whose exact form is not important for our 
analysis.  However, the dot product with $\omega$ on the far right-hand side will be important in order for us
to exploit the Coulomb gauge \eqref{CG-underline} properly; without it the argument only appears to work for sufficiently high
dimension (e.g. $n \geq 10$).

{\bf Remark.} The projection $P_k \Pi_{\omega, \theta}$, which is a smooth Fourier projection to the region $\{ \zeta \in \R^n: |\zeta| \sim 2^k, \angle(\zeta,\omega) \sim \theta \cup \angle (-\zeta,\omega)\sim \theta\}$.  On the Fourier support of this projection, the operator $\Delta_{\omega^\perp}^{-1}$ has a symbol of magnitude $O((2^k \theta)^{-2})$ on the Fourier support of $P_k \Pi_{\omega, \theta}$.  To emphasize this, we will write
$$ \Delta_{\omega^\perp}^{-1} P_k \Pi_{\omega, \theta} = (2^k \theta)^{-2} \tilde P_k \tilde \Pi_{\omega,\theta}$$
where $\tilde P_k$ and $\tilde \Pi_{\omega,\theta}$ are Fourier multipliers which obey essentially the same symbol estimates as $P_k$ and $\Pi_{\omega,\theta}$.  In particular we have
\begin{equation}\label{eq:psi-dyadic-2}
\Psi_\pm = c_\pm(\omega) \nabla_{x,t} \sum_{k < -5} \sum_{j \geq 1} (2^k  2^{\sigma k + j})^{-2} \tilde \Pi_{\omega, 2^{\sigma k+j}} \tilde P_k \underline{A}\cdot \omega.
\end{equation}

{\bf Remark.} In what follows we shall simply say that
$\angle (\zeta,\omega)\sim \theta$ on the support of the projections $\Pi_{\omega,\theta}$ and $\tilde \Pi_{\omega,\theta}$. In fact, we can redefine 
$\angle (\zeta,\omega)$ to be the minimum of the angles between $(\zeta,\omega)$
and $(-\zeta,\omega)$.

We have now decomposed $\Psi_\pm$ into a number of pieces.  To estimate these pieces
we will need to use the Coulomb gauge condition \eqref{CG-underline}.  Observe from this
gauge that, the vectorfield $F=\Box A$ is also divergence free.
We now establish an additional gain in the estimate for expressions of the form  
$\Pi_{\omega,\le \theta} B\cdot\omega $ (which certainly arise in the above decompositions)
under the condition that $B$ is divergence free.  Indeed we
gain an additional factor of $\theta$:

\begin{lemma}\label{lem:Coulomb}
Let $\omega\in {\Bbb S}^{n-1}$ be a direction. Suppose that a vectorfield $B$
 is divergence free:
 $$
 \div B = \partial_j B^j=0
 $$
 Then  for any $\theta <<1$ and an arbitrary $\zeta\in {\Bbb R}^n,\, \zeta\ne 0$,
 $$
 |\widehat {\Pi_{\omega,\le \theta} B}(\zeta)\cdot \omega|\les 
 \theta  |\widehat {\Pi_{\omega,\le \theta} B}(\zeta)|.
 $$
 In particular, for any $s\in {\Bbb R}$, 
 $$
 \|\Pi_{\omega,\le \theta} B\cdot \omega\|_{H^s_x}\les 
 \theta\, \|\Pi_{\omega,\le \theta} B\|_{H^s_x}.
 $$
 The same results hold with the projection $\Pi_{\omega,\theta}$.
\end{lemma}
\begin{proof}
On the Fourier side the Coulomb gauge  $\div B=0$ translates into the condition
$$
\hat B(\zeta) \cdot \frac {\zeta}{|\zeta|}=0,\qquad \forall \zeta\ne 0.
$$
Therefore,
\begin{align*}
\widehat {\Pi_{\omega,\le \theta} B}(\zeta)\cdot \omega& = 
\big (1-(1-\eta (\frac {\angle (\zeta,\omega)}\theta)) (1-\eta (\frac {\angle (-\zeta,\omega)}\theta)) \big)\hat B(\zeta)\cdot \omega\\ &=
\eta (\frac {\angle (\zeta,\omega)}\theta) \hat B(\zeta)\cdot \omega +
\eta (\frac {\angle (-\zeta,\omega)}\theta)) \hat B(\zeta)\cdot \omega\\ &=
\eta (\frac {\angle (\zeta,\omega)}\theta) \hat B(\zeta)\cdot 
(\omega - \frac {\zeta}{|\zeta|}) + \eta (\frac {\angle (-\zeta,\omega)}\theta) \hat B(\zeta)\cdot 
(\omega + \frac {\zeta}{|\zeta|}) 
\end{align*}
and the desired conclusion follows since $|\omega - {\zeta}/{|\zeta|}|\les \theta$
on the support of the function
$\eta (\frac {\angle (\zeta,\omega)}\theta) $ as well as 
$|\omega + {\zeta}/{|\zeta|}|\les \theta$
on the support 
$\eta (\frac {\angle (-\zeta,\omega)}\theta) $.
\end{proof}

We are now in a position to obtain a good estimate on the phase $\Psi_\pm$.

\begin{lemma}\label{lem:psid-1}  Let $\pm$ be a sign.  Then for any $2 \leq q \leq \infty$, we have
$$ \| \nabla_{x,t} \Psi_\pm(t,x,\xi) \|_{D(L^q_t L^\infty_x)} \lesssim \eps.$$
\end{lemma}

\begin{proof}
It will suffice to show that
$$ \| \nabla_{x,t}^s \Psi_\pm(t,x,\xi) \|_{D(L^2_t L^\infty_x)} \lesssim \eps$$
for $s=1,2$, since the claim for $L^q_t L^\infty_x$ then follows by a Poincare inequality type argument in time
(writing $\nabla_{x,t} \Psi_\pm$ in terms of local averages of $\nabla_{x,t} \Psi_\pm$ and $\partial_t \nabla_{x,t} \Psi_\pm$ in time, and using Minkowski's inequality).  Actually we shall just prove the claim when $s=1$; the $s=2$
case is easier due to the low frequency of $\underline{A}$ and hence of $\Psi_\pm$.

By \eqref{eq:psi-dyadic-2} and the triangle inequality, it suffices to show that
$$ \| \nabla^2_{x,t} \tilde \Pi_{\omega, 2^{\sigma k+j}} \tilde P_k \underline{A}\cdot\omega  \|_{D(L^2_t L^\infty_x)} \lesssim \eps (2^k  2^{\sigma k + j})^2 2^{\delta k}$$
for all $k < -5$; note that for fixed $k$, there are only $O(|k|)$ values of $j$ for which the summand in \eqref{eq:psi-dyadic-2} does not vanish.  Also note that the function $c_\pm(\omega)$ can be ignored thanks to \eqref{eq:factor}.
We apply Lemma \ref{banach} with $\theta := 2^{\sigma k + j}$, and reduce ourselves to showing that
$$ (\theta^{-(n-1)} \int_{S^{n-1}} \| \nabla_\xi^l \nabla^2_{x,t} \tilde \Pi_{\omega, \theta} \tilde P_k \underline{A} \cdot\omega\|_{L^2_t L^\infty_x}^2\ d\omega)^{1/2} \lesssim \eps \theta^{-l} (2^k  \theta)^2 2^{\delta k}$$
for all $0 \leq l \leq 100n$.  Note however that $(\theta\nabla_\xi)^l \tilde \Pi_{\omega,\theta} = |\xi|^{-l} \Pi^{(l)}_{\omega,\theta}$, where $\Pi^{(l)}_{\omega,\theta}$ obeys essentially the same estimates as $\Pi_{\omega,\theta}$.  The factor $|\xi|^{-l}$ can be discarded by \eqref{eq:factor}.

Fix $l$. The Fourier projection $\Pi^{(l)}_{\omega, \theta} P_k $ projects onto a region in frequency space of volume $O(2^{nk} \theta
^{n-1})$. By Bernstein's inequality, we thus have
$$ \| \nabla^2_{x,t} \Pi^{(l)}_{\omega, \theta} P_k \underline{A} \cdot\omega\|_{L^2_t L^\infty_x}
\lesssim (2^{nk} \theta^{n-1})^{\frac{1}{p_*}} \| \nabla^2_{x,t} \Pi^{(l)}_{\omega, \theta} P_k \underline{A} \cdot\omega\|_{L^2_t L^{p_*}_x}.$$
By Strichartz, this is in turn bounded by
$$ \lesssim
(2^{nk} \theta^{n-1})^{\frac{1}{p_*}} 
2^{(\frac{5}{2} - \frac{n}{p_*})k}
(\| \Pi^{(l)}_{\omega, \theta} P_k \underline{A}[0] \cdot\omega\|_{D}
+ \| \nabla_{x,t} \Pi^{(l)}_{\omega, \theta} P_k F\cdot \omega \|_{L^1_t \dot H^{n/2-3}_x}).$$
To simplify the exposition we shall just treat the contribution of the forcing term $F$, as the contribution of the initial data $\underline{A}[0]$ is similar. 
By Lemma \ref{lem:Coulomb} 
$$
\| \nabla_{x,t} \Pi^{(l)}_{\omega, \theta} P_k F\cdot \omega \|_{L^1_t \dot H^{n/2-3}_x} \les \theta\, \| \nabla_{x,t} \Pi^{(l)}_{\omega, \theta} P_k F\|_{L^1_t \dot H^{n/2-3}_x}
$$
It thus remains to show that
$$ (2^{nk} \theta^{n-1})^{\frac{1}{p_*}} 
2^{(\frac{5}{2} - \frac{n}{p_*})k}\,\theta
(\theta^{1-n} \int_{S^{n-1}} \| \nabla_{x,t} \Pi^{(l)}_{\omega, \theta} P_k F \|_{L^1_t \dot H^{n/2-3}_x}^2\ d\omega)^{1/2} \lesssim \eps (2^k  \theta)^2 2^{\delta k}.$$
By Minkowski's inequality, it suffices to show
$$ (2^{nk} \theta^{n-1})^{\frac{1}{p_*}} 
2^{(\frac{5}{2} - \frac{n}{p_*})k}\,\theta
\int (\theta^{1-n} \int_{S^{n-1}} \| \nabla_{x,t} \Pi^{(l)}_{\omega, \theta} P_k F(t) \|_{\dot H^{n/2-3}_x}^2\ d\omega)^{1/2}\ dt \lesssim \eps (2^k  \theta)^2 2^{\delta k}.$$
Observe that when two different values of $\omega$ are $\gg \theta$ apart, the corresponding functions $\nabla_{x,t} \Pi^{(l)}_{\omega, \theta} P_k F(t)$ are orthogonal.  One can then use orthogonality arguments to estimate the left-hand side by
$$ (2^{nk} \theta^{n-1})^{\frac{1}{p_*}} 
2^{(\frac{5}{2} - \frac{n}{p_*})k}\,\theta
\int \| \nabla_{x,t} P_k F(t) \|_{\dot H^{n/2-3}_x}\ dt$$
which by \eqref{eq:f-bound} is bounded by
$$ \eps (2^{nk} \theta^{n-1})^{\frac{1}{p_*}}\theta 
2^{(\frac{5}{2} - \frac{n}{p_*})k} = \eps 2^{\frac{5k}{2}} \theta^{\frac{n-1}{p_*}+1}.$$
Observe that $\frac{n-1}{p_*} \geq 1$ since $n \geq  6$ and $p_*$ is close to $2(n-1)/(n-3)$.  Since $\sigma < 1/2$,
the right-hand is thus bounded by $O(\eps (2^k \theta)^2 2^{\delta k})$ as desired.
\end{proof}

Having estimated the phase $\Psi_\pm$, we now turn to the amplitude correction $\Omega_\pm$.

\begin{lemma}\label{lem:psid-2}  Let $\pm$ be a sign.  Then the function $\Omega_\pm(t,x,\xi)$ defined in \eqref{eq:omega-def} can be decomposed as
$$ \Omega_\pm = \Omega_\pm^1 + \Omega_\pm^2,$$
where 
$$ \| \Omega_\pm^1 \|_{D(L^1_t L^\infty_x)} \lesssim \eps$$
and
$$ \| \Omega_\pm^2 \|_{D(L^2_t L^q_x)} \lesssim \eps,$$
where $q$ is the exponent such that $1/q + 1/p_* = 1/2$ (i.e. $q$ is slightly less than $n-1$).
\end{lemma}

The significance of the spaces $L^1_t L^\infty_x$ and $L^2_t L^q_x$ is that any function in those spaces, multiplied by a (frequency-localized) function in the $S$ space, will fall into the $L^1_t L^2_x$ space.

\begin{proof}
In light of \eqref{eq:omega-def} and \eqref{eq:psi-small}, it will suffice to prove the following claims:
\begin{align}
\| |\xi| \sum_{k < -5} \Pi_{\omega, \leq 2^{\sigma k}} P_k \underline{A} \cdot \omega \|_{D(L^2_t L^q_x)}  &\lesssim \eps \label{eq:d1}\\
\| |\xi| \sum_{k < -5} \Delta_{\omega^\perp}^{-1} \Pi_{\omega, > 2^{\sigma k}} P_k F \cdot \omega \|_{D(L^1_t L^\infty_x)} &\lesssim \eps  \label{eq:d2}\\
\| \Box \Psi_\pm \|_{D(L^1_t L^\infty_x)} &\lesssim \eps  \label{eq:d3}\\
\| |\nabla_x \Psi_\pm|^2 - |\partial_t \Psi_\pm|^2 \|_{D(L^1_t L^\infty_x)} &\lesssim \eps  \label{eq:d4}\\
\| \underline{A} \cdot \nabla_x \Psi_\pm \|_{D(L^1_t L^\infty_x)} &\lesssim \eps.  \label{eq:d5}
\end{align}

The quadratic estimates \eqref{eq:d4}, \eqref{eq:d5} are quite easy, thanks to the bound previously obtained
for $\Psi_\pm$.  Indeed, the estimate \eqref{eq:d4} follows immediately from Lemma \ref{lem:psid-1} (with $q=2$) and \eqref{eq:algebra}.  For \eqref{eq:d5}, observe from Strichartz' inequality and the frequency localization of $\underline{A}$ that
$$ \| \underline{A} \|_{L^2_t L^\infty_x} \lesssim \| \underline{A} \|_{L^2_t L^n_x} \lesssim \| \underline{A}[0] \|_D + \| \Box \underline{A} \|_{L^1_t \dot H^{n/2-2}_x} \lesssim \eps,$$
so the claim again follows from Lemma \ref{lem:psid-1} and \eqref{eq:algebra}.

It remains to prove the linear estimates \eqref{eq:d1}, \eqref{eq:d2}, \eqref{eq:d3}.  We begin with \eqref{eq:d1}.  We can discard the $|\xi| $ factor by \eqref{eq:factor}.  By the triangle inequality, it will suffice to show that 
$$ \| \Pi_{\omega, \leq 2^{\sigma k}} P_k \underline{A}\cdot \omega \|_{D(L^2_t L^q_x)} \lesssim \eps 2^{\delta k}$$
for all $k < -5$.  Fix $k$.  Applying Lemma \ref{banach} with $\theta := 2^{\sigma k}$, it suffices to show that
$$
(\theta^{1-n} \int_{S^{n-1}} \| (\theta\nabla_\xi)^l \Pi_{\omega, \leq \theta} P_k \underline{A}\cdot \omega 
 \|_{L^2_t L^q_x}^2 \ d\omega)^{1/2} \lesssim \eps 2^{\delta k}$$
for all $0 \leq l \leq 100n$.  

Fix $l$.  As before, we write $(\theta\nabla_\xi)^l \Pi_{\omega, \leq \theta} = |\xi|^{-l} \Pi^{(l)}_{\omega, \leq \theta}$ and discard the $|\xi|^{-l}$ factor.  Thus it will suffice to show that
\begin{equation}\label{2q}
(\theta^{1-n} \int_{S^{n-1}} \| \Pi^{(l)}_{\omega, \leq \theta} P_k \underline{A}\cdot \omega
 \|_{L^2_t L^q_x}^2\ d\omega)^{1/2} \lesssim \eps 2^{\delta k}.
\end{equation}
The Fourier projection $\Pi^{(l)}_{\omega, \leq 2^{\sigma k}} P_k $ projects onto a region in frequency space of volume $O(2^{nk} \theta^{n-1})$. By Bernstein's inequality, we thus have
$$ \| \Pi^{(l)}_{\omega, \leq \theta} P_k \underline{A}\cdot \omega
 \|_{L^2_t L^q_x} \lesssim (2^{nk} \theta^{n-1})^{\frac{1}{p_*} - \frac{1}{q}} \| \Pi^{(l)}_{\omega, \leq \theta} P_k \underline{A}\cdot \omega \|_{L^2_t L^{p_*}_x};$$
note that $p_* \le q$ since $n \geq 6$ and $\delta$ is small.  By Strichartz' inequality, we have
$$ \| \Pi^{(l)}_{\omega, \leq \theta} P_k \underline{A}\cdot\omega \|_{L^2_t L^{p_*}_x}
\lesssim 2^{(\frac{1}{2} - \frac{n}{p_*}) k} 
(\| \Pi^{(l)}_{\omega, \leq \theta} P_k \underline{A}[0]\cdot\omega \|_{D}
+ \| \Pi^{(l)}_{\omega, \leq \theta} P_k F\cdot\omega \|_{L^1_t \dot H^{n/2-2}_x}).$$
Again we shall just treat the contribution of the forcing term $F$.  We can now bound the left-hand side of \eqref{2q}, by
$$
\lesssim (2^{nk} \theta^{n-1})^{\frac{1}{p_*} - \frac{1}{q}} 2^{(\frac{1}{2} - \frac{n}{p_*}) k} 
(\theta^{1-n} \int_{S^{n-1}} \| \Pi^{(l)}_{\omega, \leq \theta} P_k F\cdot\omega \|_{L^1_t \dot H^{n/2-2}_x}^2\ d\omega)^{1/2}.$$
By Lemma \ref{lem:Coulomb} and Minkowski's inequality we may bound this by
$$
\lesssim (2^{nk} \theta^{n-1})^{\frac{1}{p_*} - \frac{1}{q}} 2^{(\frac{1}{2} - \frac{n}{p_*}) k} 
\,\theta \int (\theta^{1-n} \int_{S^{n-1}} \| \Pi^{(l)}_{\omega, \leq \theta} P_k F(t)  \|_{\dot H^{n/2-2}_x}^2\ d\omega)^{1/2}\ dt.
$$
Observe that when two different values of $\omega$ are $\gg \theta$ apart, the corresponding functions $\Pi^{(l)}_{\omega, \leq \theta} P_k F(t)$ are orthogonal.  One can then use orthogonality arguments to estimate the previous expression by
$$
\lesssim (2^{nk} \theta^{n-1})^{\frac{1}{p_*} - \frac{1}{q}}2^{(\frac{1}{2} - \frac{n}{p_*}) k} \,
\theta
\int \| P_k F(t) \|_{\dot H^{n/2-2}_x}\ dt.
$$
But by \eqref{eq:f-bound} we have $\int \| P_k F(t) \|_{\dot H^{n/2-2}_x}\ dt \lesssim \eps$.  Thus, substituting $\theta = 2^{\sigma k}$, we will be done as long as
$$ (n + (n-1)\sigma) (\frac{1}{p_*} - \frac{1}{q}) + (\frac{1}{2} - \frac{n}{p_*}) +\sigma  > \delta.$$
If $\delta$ is sufficiently small, this will be achieved as long as
$$ (n + (n-1)\sigma) (\frac{n-3}{2(n-1)} - \frac{1}{n-1}) + (\frac{1}{2} - \frac{n(n-3)}{2(n-1)})
+\sigma  > 0.$$
After some algebra, this becomes equivalent to
$$ \frac{n-3}{2} \sigma > \frac{n+1}{2(n-1)},$$
which follows from \eqref{sigma-bounds}.  This proves \eqref{eq:d1}.

Now we prove \eqref{eq:d2}.  We decompose the left-hand side as
$$ \| |\xi| \sum_{k < -5} \sum_{j \geq 1} (2^k  2^{\sigma k + j})^{-2} \tilde \Pi_{\omega, 2^{\sigma k+j}} \tilde P_k F \cdot \omega\|_{D(L^1_t L^\infty_x)}$$
(cf. \eqref{eq:psi-dyadic-2}); of course the summand is only non-zero when $2^{\sigma k + j} \lesssim 1$.  By \eqref{eq:factor} we can discard the $|\xi|$ factor.  By \eqref{eq:f-bound}  and the triangle inequality, it will suffice to show that 
$$ \| (2^k  2^{\sigma k + j})^{-2}  \tilde \Pi_{\omega, 2^{\sigma k+j}} \tilde P_k F \cdot\omega\|_{D(L^1_t L^\infty_x)}
\lesssim (2^{\sigma k + j})^\delta \| \tilde P_k F \|_{L^1_t \dot H^{n/2-2}_x}$$
for each $k < -5$ and $j \geq 1$ (note how it is vital here that $F$ was controlled in an $l^1$ Besov space instead of an $l^2$ Besov space).  

Fix $j,k$.  Applying Lemma \ref{banach} with $\theta := 2^{\sigma k + j}$, it suffices to show that
$$
(\theta^{1-n} \int_{S^{n-1}} \| (\theta \nabla_\xi)^l \tilde \Pi_{\omega, \theta} \tilde P_k F \cdot\omega \|_{L^1_t L^\infty_x}^2\ d\omega)^{1/2} \lesssim (2^k \theta)^2 \theta^\delta  \| \tilde P_k F \|_{L^1_t \dot H^{n/2-2}_x}$$
for all $0 \leq l \leq 100n$.  

Fix $l$.  Similarly to before, we may write $(\theta \nabla_\xi)^l \tilde \Pi_{\omega, \theta} = |\xi|^{-l} \tilde \Pi^{(l)}_{\omega,\theta}$, where $\tilde \Pi^{(l)}_{\omega,\theta}$ obeys similar estimates to $\tilde \Pi_{\omega,\theta}$.  It thus suffices to show that
$$
(\theta^{1-n} \int_{S^{n-1}} \| \tilde \Pi^{(l)}_{\omega, \theta} \tilde P_k F \cdot\omega \|_{L^1_t L^\infty_x}^2\ d\omega)^{1/2} \lesssim (2^k \theta)^2 \theta^\delta \| \tilde P_k F \|_{L^1_t \dot H^{n/2-2}_x}.$$
By Minkowski it suffices to show
\begin{equation}\label{theta-stuff}
\int (\theta^{1-n} \int_{S^{n-1}} \| \tilde \Pi^{(l)}_{\omega, \theta} \tilde P_k F(t) \cdot\omega \|_{L^\infty_x}^2\ d\omega)^{1/2}\ dt \lesssim (2^k \theta)^2 \theta^\delta \| \tilde P_k F \|_{L^1_t \dot H^{n/2-2}_x}.
\end{equation}
By Bernstein and Lemma \ref{lem:Coulomb} we have
$$ \| \tilde \Pi^{(l)}_{\omega, \theta} \tilde P_k F(t)\cdot\omega  \|_{L^\infty_x} \lesssim 2^{nk/2} \theta^{(n-1)/2} \| \tilde \Pi^{(l)}_{\omega, \theta} P_k F(t)\cdot\omega  \|_{L^2_x}
\sim 2^{2k} \theta^{(n+1)/2} \| \tilde \Pi^{(l)}_{\omega, \theta} P_k F(t) \|_{\dot H^{n/2 - 2}_x},$$
so by orthogonality we can bound the left-hand side of \eqref{theta-stuff} by
$$
2^{2k} \theta^{(n+1)/2} \| \tilde P_k F(t) \|_{\dot H^{n/2 - 2}_x},$$
which is acceptable since $n \geq \n$ and $\theta$ is small.  This proves \eqref{eq:d2}.

Now we prove \eqref{eq:d3}.  By \eqref{eq:psi-dyadic-2} we can write this expression as
$$
\Box \Psi_\pm = c_\pm(\omega) \nabla_{x,t} \sum_{k < -5} \sum_{j \geq 1} (2^k  2^{\sigma k + j})^{-2} \tilde \Pi_{\omega, 2^{\sigma k+j}} \tilde P_k F\cdot\omega.$$
By the triangle inequality and \eqref{eq:factor}, it will thus suffice to show that 
$$ \| \nabla_{x,t} \tilde \Pi_{\omega, 2^{\sigma k+j}} \tilde P_k F\cdot\omega \|_{D(L^1_t L^\infty_x)} \lesssim 2^{(2+\delta)k} 2^{2(\sigma k+j)}\epsilon$$
for all $k < -5$ and $j \geq 1$, 
since the claim will then follow by summing in $k$ and $j$ (note that for each $k$ there are only $O(k)$ values of $j$ for which the summand is non-zero).  Applying Lemma \ref{banach} with $\theta := 2^{\sigma k + j}$ and arguing as before, it suffices to show that
$$
(\theta^{1-n} \int_{S^{n-1}} \| \nabla_{x,t} \tilde \Pi^{(l)}_{\omega, \theta} \tilde P_k F
\cdot\omega \|_{L^1_t L^\infty_x}^2\ d\omega)^{1/2} \lesssim 2^{(2+\delta)k} \theta^2\epsilon.$$
Again using Minkowski, it suffices to show
$$
\int (\theta^{1-n} \int_{S^{n-1}} \| \nabla_{x,t} \tilde \Pi_{\omega, \theta} \tilde P_k F(t) \cdot\omega \|_{L^\infty_x}^2\ d\omega)^{1/2}\ dt \lesssim 2^{(2+\delta)k} \theta^2 \epsilon.$$
By Bernstein and Lemma \ref{lem:Coulomb} again, we have
\begin{align*}
 \| \nabla_{x,t} \tilde \Pi_{\omega, \theta} \tilde P_k F(t) \cdot\omega\|_{L^\infty_x} &\lesssim \theta^{(n-1)/2}
2^{3k} \| \nabla_{x,t} \tilde \Pi_{\omega, \theta} \tilde P_k F(t)\cdot\omega \|_{\dot H^{n/2 - 3}_x}\\ &\lesssim \theta^{(n+1)/2}
2^{3k} \| \nabla_{x,t} \tilde \Pi_{\omega, \theta} \tilde P_k F(t) \|_{\dot H^{n/2 - 3}_x}
\end{align*}
and so by orthogonality it suffices to show
$$
2^{3k} \theta^{(n+1)/2} \| \tilde P_k \nabla_{x,t} F \|_{L^1_t \dot H^{n/2-3}_x} \lesssim 2^{(2+\delta)k} \theta^2 \epsilon.$$
But this follows from \eqref{eq:f-bound}.  This proves \eqref{eq:d3}.

We have now estimated all the terms which appear in $\Omega_\pm$, and the proof of the Lemma is complete.

\end{proof}

\section{Approximate unitarity of $U^\pm(t)$}\label{sec:unitarity}

In this section, the main property of $\underline{A}$ (apart from the Coulomb gauge condition
\eqref{CG-underline}) is simply the Sobolev bound
\begin{equation}\label{eq:a-sob}
 \| \nabla^2_{x,t} \underline{A}(t) \|_{\dot H^{n/2-3}_x} \lesssim \eps
\end{equation}
for all times $t \in \R$.  To see this, first observe from \eqref{eq:f-sup} that
$$
 \|\Box \underline{A}(t) \|_{\dot H^{n/2-3}_x} \lesssim \eps
$$
so it suffices to show that
$$ 
\| \nabla_x \nabla_{x,t} \underline{A}(t) \|_{\dot H^{n/2-3}_x} \lesssim \eps.$$
But this follows from \eqref{a-small} and energy estimates.

The goal of this section is to prove $L^2$ bounds on the operators $U^\pm(t)$.  

\begin{proposition}\label{prop:l2}  For any $t \in I$ and sign $\pm$, we have
$$ \| U^\pm(t) h\|_{L^2_x} \lesssim \| h \|_2$$
\end{proposition}

\begin{proof}
Fix the time $t$ and the sign $\pm$; we will allow all the quantities we define to depend implicitly on $t$ and
$\pm$.  By the $TT^*$ method, it will suffice to show the estimate
$$ \| U^\pm(t) U^\pm(t)^* f \|_{L^2_x} \lesssim \| f \|_2,$$
where $U^\pm(t)^*$ is the adjoint of $U^\pm(t)$ with respect to the standard inner products in the $x$ and $\xi$ variables.
A direct calculation shows that
\begin{equation}\label{eq:uu-form}
 U^\pm(t) U^\pm(t)^* f(x) := \int K(x,y) f(y)\ dy
\end{equation}
where $K(x,y)$ is the kernel
$$
 K(x,y) := \int e^{2\pi i \Phi(x,y,\xi)} e^{2\pi i (x-y) \cdot \xi} 
a^2(\xi)\ d\xi
$$
and $\Phi$ is the phase
$$ \Phi(x,y,\xi) := \Psi_\pm(t,x,\xi) - \Psi_\pm(t,y,\xi)$$
(we are suppressing the dependence on $t$ and $\pm$ as these quantities are fixed).
By Schur's test and the self-adjointness of the kernel $K$, it will thus suffice to show the  estimate
\begin{equation}\label{schur}
\int |K(x,y)|\, dy\les 1 \hbox{ for all } x \in \R^n.
\end{equation}
Observe that $K$ is trivially bounded by replacing everything by absolute values, so we may restrict the integral
to the region where $|x-y| \gg 1$.

We now estimate the kernel $K(x,y)$ on the region $|x-y| \gg 1$.  We write $x-y$ in polar co-ordinates as
$$
(x-y) = |x-y|\, \alpha,\qquad \alpha = \frac {x-y}{|x-y|}
$$
Let $b(\lambda)$  be a smooth cut-off function supported on $[-2,2]$ which equals 1 on $[-1,1]$.  Define angle $\theta_*$ 
 $$\theta_* := |x-y|^{-1+2\delta}.
 $$ 
 Then
we can split 
\begin{align}
K(x,y) =
& \int e^{2\pi i\Phi (x,y,\xi)} e^{2\pi i (x-y) \cdot \xi} a^2(\xi)\,b( \frac{\xi\c\alpha}{\sqrt{\theta_*}})
 d\xi + \label{K-1}\\
 &\int e^{2\pi i\Phi (x,y,\xi)} e^{2\pi i (x-y) \cdot \xi} a^2(\xi)\,(1-b(\frac{\xi\c\alpha}{\sqrt{\theta_*}} ))
 d\xi \label{K-2}
\end{align}
We consider first the term \eqref{K-2}, where the angle between $\xi$ and $x-y$ is less than $\frac{\pi}{2} -O(\sqrt{\theta_*})$.
This term is relatively easy to handle and one only needs to exploit the fact that the phases $\Phi_\pm (t,x,\xi)$ 
are homogeneous functions in $|\xi|$ of degree $0$.
On the support of $a^2(\xi)\,(1-b(\frac{\xi\c\alpha}{\sqrt{\theta_*}} ))$ we have 
$$
\big | (x-y) \cdot \frac {\xi}{|\xi|}\big |\ge \frac 12 |x-y|^{\frac 12+\delta}
$$
while 
$$
|\nabla^K_{|\xi|} \big (a^2(\xi) (1-b( \frac{\xi\c\alpha}{\sqrt{\theta_*}} )\big )|\les
|x-y|^{K(\frac 12-\delta)}
$$
for all $K\ge 0$.
Therefore,  integrating by parts with respect to $|\xi|$ we derive that
\begin{equation}\label{k-bound}
|\eqref{K-2}| 
\les |x-y|^{-2\de K} \les |x-y|^{-n-\sqrt\delta}
\end{equation}
for a sufficiently large $K>0$.  Integrating this in $y$ we thus see that the contribution of this term to
\eqref{schur} is acceptable.

It remains to deal with the main term \eqref{K-1}, when 
$\xi$ is localized to make an angle of $\frac{\pi}{2} + O(\sqrt{\theta_*})$ with
$x-y$.  In this case we shall in fact obtain the bound
\begin{equation}\label{k-bound-2}
 \int_{|x-y| \gg 1} |\eqref{K-1}|\ dy \lesssim \eps.
\end{equation}

To prove \eqref{k-bound-2} we use \eqref{eq:psi-dyadic} to decompose the phase $\Phi(t,x,y,\xi)$ as
\begin{equation}\label{eq:Phi-decomp}
\Phi(x,y,\xi) = \Phi_{smooth}(x,y,\xi) + \Phi_{small}(x,y,\xi), 
\end{equation}
where $\Phi_{smooth}$ and $\Phi_{small}$ are defined as
$$
\Phi_{small}(x,y,\xi) := \Psi_\pm^{<\theta_*}(t,x,\xi) - \Psi_\pm^{<\theta_*}(t,y,\xi),
\qquad 
\Phi_{smooth}(x,y,\xi) := \Psi_\pm^{\geq \theta_*}(t,x,\xi) - \Psi_\pm^{\geq \theta_*}(t,y,\xi)
$$
where $\theta_*$ is the angle $\theta_* := |x-y|^{-1+2\delta}$ and $\Psi_\pm^{<\theta_*}$ is the function
\begin{equation}\label{psi-less}
 \Psi_\pm^{<\theta_*} :=
c^\pm(\omega) \cdot \nabla_{x,t} \sum_{k < -5} 
\sum_{j \geq 1: 2^{\sigma k + j} < \theta_*} 
(2^k  2^{\sigma k + j})^{-2} \tilde \Pi_{\omega, 2^{\sigma k + j}} \tilde P_k \underline{A} 
\cdot \omega
\end{equation}
and similarly
\begin{equation}\label{psi-more}
 \Psi_\pm^{\geq \theta_*} :=
c^\pm(\omega) \cdot \nabla_{x,t} \sum_{k < -5} 
\sum_{j \geq 1: 2^{\sigma k + j} \geq \theta_*} 
(2^k  2^{\sigma k + j})^{-2} \tilde \Pi_{\omega, 2^{\sigma k + j}} \tilde P_k \underline{A} 
\cdot \omega.
\end{equation}

From \eqref{eq:Phi-decomp} we have in particular that
$$ e^{2\pi i \Phi_\pm(x,y,\xi)} = e^{2\pi i \Phi_{smooth}(x,y,\xi)} + O(|\Phi_{small}(x,y,\xi)|),$$
so to complete the proof of \eqref{schur} (and hence the Proposition) it will suffice to show the bound
\begin{equation}\label{eq:Phi-small-bound}
\int_{|x-y| \gg 1} (\int |\Phi_{small} (x,y,\xi)| a^2(\xi) \, b(\frac{\xi\c\alpha}{\sqrt{\theta_*}}) d\xi )\,dy
\les \epsilon
\end{equation}
for all $x \in \R^n$, and the bound
\begin{equation}\label{eq:Phi-smooth-bound}
|\int e^{2\pi i\Phi_{smooth}(x,y,\xi)} e^{2\pi i (x-y) \cdot \xi} a^2(\xi)\,
b(\frac{\xi \cdot \alpha}{\sqrt{\theta_*}})  d\xi  |\les |x-y|^{-n-\sqrt\delta}
\end{equation}
for all $x,y \in \R^n$ with $|x-y| \gg 1$.

\divider{Proof of \eqref{eq:Phi-small-bound}.}

We first show \eqref{eq:Phi-small-bound}.  Fix $x$.  By the mean-value theorem we have
$$ |\Phi_{small}(x,y,\xi)| \lesssim |x-y| \| (\alpha \cdot \nabla_{x'}) \Psi_\pm^{< \theta_*}(t) \|_{L^\infty_{x'}}$$
where we have named the spatial variable $x'$ in order to distinguish it from the point $x$.  By
\eqref{psi-less} and the triangle inequality we thus have
$$ |\Phi_{small}(x,y,\xi)| \lesssim |x-y|
\sum_{k < -5} \sum_{j \geq 1: 2^{\sigma k + j} < \theta_*} (2^k  2^{\sigma k + j})^{-2}
\| (\alpha \cdot \nabla_{x'}) \nabla_{x',t}
 \tilde \Pi_{\omega, 2^{\sigma k + j}} \tilde P_k \underline{A} 
\cdot \omega \|_{L^\infty_{x'}}.$$
Estimating $\| f\|_{L^\infty_{x'}} \lesssim \| \hat f\|_{L^1_\xi}$, we thus have
$$ |\Phi_{small}(x,y,\xi)| \lesssim |x-y|
\sum_{k < -5} \sum_{j \geq 1: 2^{\sigma k + j} < \theta_*} (2^k  2^{\sigma k + j})^{-2}
\int_{|\zeta| \sim 2^k; \angle(\zeta,\omega) \sim 2^{\sigma k + j} } |((\alpha\c\nabla_{x'}) \nabla_{x',t} \underline{A})\hat{} \cdot \omega(t,\zeta)|\ d\zeta.$$
The expression $(\alpha\c\nabla_{x'}) \nabla_{x',t} \underline{A}$ is divergence-free in $x'$, so we can
use Lemma \ref{lem:Coulomb} to obtain
$$ |\Phi_{small}(x,y,\xi)| \lesssim |x-y|
\sum_{k < -5} \sum_{j \geq 1: 2^{\sigma k + j} < \theta_*} 2^{-2k}  2^{-\sigma k - j}
\int_{|\zeta| \sim 2^k; \angle(\zeta,\omega) \sim 2^{\sigma k + j} } |((\alpha\c\nabla_{x'}) \nabla_{x',t} \underline{A})\hat{}|\ d\zeta.$$
Thus we have
\begin{align*}
\int |\Phi_{small} (\xi)| a^2(\xi) \, &b( \frac{\xi\c\alpha}{\sqrt{\theta_*}}) d\xi\les
|x-y| 
\sum_{k < -5} 
\sum_{j \geq 1: 2^{\sigma k + j} < \theta_*} 2^{-2k} 2^{-(\sigma k + j)}\c\\&\c
\int_{\angle(\omega,\alpha)=\frac{\pi}{2} + O(\sqrt{\theta_*})}
\int_{|\zeta| \sim 2^k; \angle(\zeta,\omega) \sim 2^{\sigma k + j} } |\widehat{\alpha\c\nabla_x \nabla_{x,t} \underline{A}}(t,\zeta)|\ d\zeta
d\omega.
\end{align*}
We have 
$$
|\widehat{\alpha\c\nabla_x \nabla_{x,t} \underline{A}}(t,\zeta)|=
|\alpha\c\zeta|\, |\widehat{\nabla_{x,t} \underline{A}}(t,\zeta)|\les
|\frac{\pi}{2}-\angle(\alpha,\zeta)|\, |\widehat{\nabla^2_{x,t} \underline{A}}(t,\zeta)|
$$
Since $\angle(\omega,\alpha)=\frac{\pi}{2} + O(\sqrt{\theta_*})$ and 
either $\angle(\zeta,\omega) \sim 2^{\sigma k + j}\le \theta_*$ or 
$\angle(-\zeta,\omega) \sim 2^{\sigma k + j}\le \theta_*$, we obtain 
that $\angle(\alpha,\zeta)=\frac{\pi}{2} + O(\sqrt{\theta_*})$ .
Interchanging the integrals, and then performing the $\omega$ integral, we obtain the bound
\begin{align}
\int |\Phi_{small} (\xi)| a^2(\xi) \, &b(\frac{\xi\c\alpha}{\sqrt{\theta_*}}) d\xi\les
|x-y| 
\sum_{k < -5} 
\sum_{j \geq 1: 2^{\sigma k + j} < \theta_*} 2^{-2k} {2^{(\sigma k + j)(n-2)}}\, 
\sqrt{\theta_*}\c\nn\\
&\c 
\int_{|\zeta| \sim 2^k; \angle(\alpha,\zeta)=\frac{\pi}{2} + O(\sqrt{\theta_*})} |\widehat{\nabla^2_{x,t} \underline{A}}(t,\zeta)|\ d\zeta.
\label{eq:Phi-small-bound-alt}
\end{align}
Integrating \eqref{eq:Phi-small-bound-alt} in $y$ and interchanging the limits
of integrations we derive
\begin{align*}
\int (\int |\Phi_{small} (\xi)| a^2(\xi) \, &b(\frac{\xi\c\alpha}{\sqrt{\theta_*}}) d\xi)\,dy\les
\int_{|\zeta| \sim 2^k}  (\int_{\angle(\alpha,\zeta)=\frac{\pi}{2} + O(\sqrt{\theta_*})}  |x-y|\c\\ &\sum_{k < -5}  
\sum_{j \geq 1: 2^{\sigma k + j} < \theta_*} 2^{-2k} {2^{(\sigma k + j)(n-2)}}\, 
|\widehat{\nabla^2_{x,t} \underline{A}}(t,\zeta)| dy ) d\zeta
\end{align*}
We now recall that $\alpha = (x-y)/|x-y|\in S^{n-1}$ and introduce on $\R^n_y$ the 
standard polar coordinates $(r,a,\beta)$, where $r=|x-y|$,
$a\in [-\pi,\pi]$, and $\beta\in S^{n-2}_{|\sin a|}$. We choose $a$ to be the angle 
between the vectors $\alpha$ and $\zeta$. 
Therefore,
$$
\int_{\angle(\alpha,\zeta)=\frac{\pi}{2} + O(\sqrt{\theta_*})}  \sin a da 
d\omega_{n-2}\les \sqrt{\theta_*}.
$$
Thus we derive
\begin{align*}
\int(\int |\Phi_{small} (\xi)| a^2(\xi) \, &b(\frac{\xi\c\alpha}{\sqrt{\theta_*}}) d\xi)\,dy\les
\int_1^\infty r\sum_{k < -5}  
\sum_{j \geq 1: 2^{\sigma k + j} < \theta_*} 2^{-2k} {2^{(\sigma k + j)(n-2)}}\, 
\c\\&\c 
 \int_{|\zeta| \sim 2^k}  \sqrt{\theta_*}
|\widehat{\nabla^2_{x,t} \underline{A}}(t,\zeta)|\ d\zeta\,  \,r^{n-1}\,dr
\end{align*}
Using  Plancherel, Cauchy-Schwarz, and \eqref{eq:a-sob} we have that
\begin{align*}
 \int_{|\zeta| \sim 2^k}  \sqrt{\theta_*}|\widehat{\nabla^2_{x,t} \underline{A}}(t,\zeta)|\ d\zeta\les 
{\sqrt{\theta_*}} 2^{kn/2}
\big ( \int_{|\zeta| \sim 2^k} |\widehat{\nabla^2_{x,t} \underline{A}}(t,\zeta)|^2
\, d\zeta\big )^{\frac 12}
\les \epsilon \,2^{3k}  {\sqrt{\theta_*}} 
\end{align*}
Thus
\begin{align*}
\int(\int |\Phi_{small} (\xi)| a^2(\xi) \, b(\sqrt{\theta_*} \xi\c\alpha) d\xi)\,dy
&\les \epsilon\, \int_1^\infty r^n\sum_{k < -5} 
\sum_{j \geq 1: 2^{\sigma k + j} < \theta_*} 2^{3k} 2^{-2k} 2^{(\sigma k + j)(n-2)}
\theta_* \, dr
\\  & \lesssim \epsilon\,\int_1^\infty r^n
\sum_{k: 2^{\sigma k} < \theta_*}  2^k \theta_*^{n-1}\,dr\\ &\les
\epsilon \int_1^\infty r^{n-(n-1+1/\sigma)(1-2\de)}\, \, dr\les \epsilon,
\end{align*}
provided that 
$n-(n-1+1/\sigma) <-1,$
which follows from the condition $\si<1/2$.  This completes the proof of \eqref{eq:Phi-small-bound}.

\divider{Proof of \eqref{eq:Phi-smooth-bound}}

We now prove \eqref{eq:Phi-smooth-bound}.  It will suffice to prove the smoothness bounds
\begin{equation}\label{smooth-est}
|\nabla_\xi^K \Phi_{smooth}(x,y,\xi)|\les |x-y|^{K(1-\de)}
\end{equation}
for all $|x-y| \gg 1$ and all
$\xi$ in the support of $a$, since this will imply (by repeated applications of the chain
and Leibnitz rules)
$$
|\nabla_\xi^K  \big (a^2(\xi)\,b( \frac{\xi\c\alpha}{\sqrt{\theta_*}})\, 
e^{2\pi i\Phi_{smooth}(\xi)}\big )|\les |x-y|^{K(1-\de)},
$$
and the bound \eqref{eq:Phi-smooth-bound} then follows with the help
of multiple integration by parts with respect to $\xi$.

It remains to prove \eqref{smooth-est}.  Fix $x,y,\xi$.  Observe that we have the two estimates
$$ |f(x) - f(y)| \lesssim \min( |x-y| \| \nabla_{x'} f \|_{L^\infty_{x'}}, \| f \|_{L^\infty_{x'}} ),$$
where we again use $x'$ to denote the spatial variable in order to distinguish it from $x$ or $y$.  In
other words, we have
$$ |f(x) - f(y)| \lesssim \min_{m=0,1} |x-y|^m \| \nabla^m_{x'} f \|_{L^\infty_{x'}}.$$
By \eqref{psi-more} and the triangle inequality, we can thus estimate the left-hand side of \eqref{smooth-est} as
$$
\lesssim \sum_{k < -5} \sum_{j \geq 1: 2^{\sigma k + j} \geq \theta_*} 
\min_{m=0,1} \| \nabla_\xi^K |x-y|^m \nabla^m_{x'} (c_\pm(\omega) \cdot \nabla_{x',t} (2^k 2^{\sigma k + j})^{-2} \tilde \Pi_{\omega, 2^{\sigma k + j}} \tilde P_k \underline{A}) \cdot\omega\|_{L^\infty_{x'}}.$$
We use the Fourier inversion formula, and observe that each derivative in $\xi$ increases the symbols of the multipliers by at most $O(2^{-\sigma k - j})$, primarily because of the projection $\tilde \Pi_{\omega, 2^{\sigma k + j}}$.  Thus the effect of $\nabla_\xi^K$ is basically to multiply things by $2^{-K(\sigma k + j)}$.  There is the
danger that one or more of the $\nabla_\xi$ derivatives hits the dot product with $\omega$, which means
that we can no longer use Lemma \ref{lem:Coulomb} to obtain the additional factor of $2^{\sigma k + j}$, but
in that case the loss arising from $\nabla_\xi^K$ is at most $2^{-(K-1)(\sigma k + j)}$ because only $K-1$ of
those derivatives can hit the projection $\tilde \Pi_{\omega, 2^{\sigma k+j}}$.  Thus in either case the combined
effect of the $\nabla_\xi^K$ derivative and the dot product with $\omega$ is to add a factor of $2^{-K(\sigma k+j)}
2^{\sigma k+j}$. Thus, by arguing as in the proof of \eqref{eq:Phi-small-bound}, we can bound this expression by
$$
\lesssim \sum_{k < -5} \sum_{j \geq 1: 2^{\sigma k + j} \geq \theta_*} 
\min_{m=0,1} \int_{|\zeta| \sim 2^k; \angle(\zeta,\omega) \sim 2^{\sigma k + j}}
(2^{-\sigma k - j})^K |x-y|^m 2^{km} 2^{-2k} 2^{-(\sigma k + j)}
|\widehat{\nabla_{x',t} \underline{A}}(\zeta)|\ d\zeta.$$
We may of course restrict the $j$ summation to the region where $2^{\sigma k + j} \lesssim 1$, since the integral is vacuous otherwise.

By Cauchy-Schwarz, Plancherel, and \eqref{eq:a-sob} we have that
\begin{equation}\label{repeat}
\begin{split}
 \int_{|\zeta| \sim 2^k; \angle(\zeta,\omega) \sim 2^{\sigma k + j}}
|\widehat{\nabla_{x',t} \underline{A}}(\zeta)|\ d\zeta 
&\lesssim (2^{(n-1)(\sigma k + j)} 2^{nk})^{1/2} (\int_{|\zeta| \sim 2^k}
|\widehat{\nabla_{x',t} \underline{A}}(\zeta)|^2\ d\zeta)^{1/2} \\
&\lesssim 2^{\frac{n-1}{2}(\sigma k + j)} 2^{2k} \| \nabla_{x,t} \underline{A} \|_{\dot H^{n/2-2}_x} \\
&\lesssim 2^{\frac{n-1}{2}(\sigma k + j)} 2^{2k} \epsilon.
\end{split}
\end{equation}
Thus we can bound the previous expression by
\begin{equation}\label{repeat-2}
\lesssim \sum_{k < -5} \sum_{j \geq 1: \theta_* \leq 2^{\sigma k + j} \lesssim 1} 
\min_{m=0,1} (2^{-\sigma k - j})^K |x-y|^m 2^{km} 2^{-2k} 2^{-(\sigma k + j)}
2^{\frac{n-1}{2}(\sigma k + j)} 2^{2k} \epsilon.
\end{equation}
We estimate $(2^{-\sigma k - j})^K$ by $\theta_*^{-K}$, so that the previous is bounded by
$$
\lesssim \epsilon \theta_*^{-K} \sum_{k < -5} \min(1, 2^k |x-y|) \sum_{j: 2^{\sigma k + j} \lesssim 1} 
2^{\frac{n-3}{2}(\sigma k + j)}.$$
Summing in $j$, we can bound this by
$$
\lesssim \epsilon \theta_*^{-K} \sum_{k < -5} \min(1, 2^k |x-y|) \lesssim
\epsilon \theta_*^{-K} \log(|x-y|),$$
which is acceptable by definition of $\theta_*$.  This proves \eqref{eq:Phi-smooth-bound}, which completes
the derivation of \eqref{schur}, and the Proposition is proved.
\end{proof}

Next, we modify the above argument to accommodate spacetime derivatives.  Fortunately, the $D(L^q_t L^r_x)$ type
control we established earlier will make this very easy.

\begin{proposition}\label{prop:l2-deriv}  For any $t \in I$ and any sign $\pm$, we have
$$ \| \nabla_x U^\pm(t) h - U^\pm(t) 2\pi i \xi h \|_{L^2_x} \lesssim \eps \| h \|_2$$
and
\begin{equation}\label{eq:t-small}
 \| \partial_t U^\pm(t) h \mp U^\pm(t) 2\pi i |\xi| h \|_{L^2_x} \lesssim \eps \| h \|_2.
\end{equation}
In particular, from these estimates and Proposition \ref{prop:l2} we have
$$ \| \nabla_{x,t} U^\pm(t) h \|_{L^2_x} \lesssim \| h \|_2$$
(note that the multiplier $\xi$ is bounded on the support of $a(\xi)$).
\end{proposition}

\begin{proof}
Fix $t$ and $\pm$.  From \eqref{eq:ut-def} we see that it suffices to show that
$$ \| \int (\nabla_{x,t} \Psi_\pm(t,x,\xi)) e^{2\pi i\Psi_\pm(t,x,\xi)} e^{2\pi i x \cdot \xi} e^{\pm 2\pi i t|\xi|} h(\xi) a(\xi)\ d\xi \|_{L^2_x} \lesssim \eps \| h \|_2.$$
However, from Proposition \ref{prop:l2} we already have
$$ \| \int e^{2\pi i\Psi_\pm(t,x,\xi)} e^{2\pi i x \cdot \xi} e^{\pm 2\pi i t|\xi|} h(\xi) a(\xi)\ d\xi \|_{L^2_x} \lesssim \| h \|_2.$$
Meanwhile, from Lemma \ref{lem:psid-1}, we have
$$ \| \nabla_{x,t} \Psi_\pm(t,x,\xi) \|_{D(L^\infty_t L^\infty_x)} \lesssim \eps.$$
The claim follows\footnote{Strictly speaking, the definition of the $D(L^\infty_t L^\infty_x)$ norm requires that we use spacetime norms rather than spatial norms on the left-hand side of our estimates.  However, this can be fixed by the artificial expedient of replacing the $L^2_x$ norm by the $L^1_t L^2_x$ norm and placing a Dirac mass at the time $t$ inside the norm.}.  
\end{proof}

We now apply the above machinery to obtain the estimates we need for \eqref{eq:data-approx}:

\begin{proposition}\label{prop:unitary}  Let $f, g$ be functions on $\R^n$ with Fourier support in $\{ 2^{-3} \leq |\xi| \leq 2^3 \}$ obeying \eqref{eq:fg-norm}.  Then there exist functions $h_+$, $h_-$ on $\R^n$ such that
\begin{equation}\label{eq:l2}
 \| h_+ \|_2 + \| h_- \|_2 \lesssim 1
\end{equation}
and
$$ \| (U^+(0) h_+ + U^-(0) h_-) - f \|_2 \lesssim \eps^{\delta}$$
and
$$ \| (\dot U^+(0) h_+ + \dot U^-(0) h_-) - g \|_2 \lesssim \eps^{\delta}$$
where $\dot U^\pm(t) := \partial_t U^\pm(t)$.
\end{proposition}

\begin{proof}
We shall construct $h_+$ and $h_-$ explicitly, as
$$ h_\pm(\xi) := \frac{1}{2} ( U^\pm(0)^* f(\xi) \pm \frac{1}{2\pi i |\xi|} U^\pm(0)^* g(\xi) ).$$
The $L^2$ bounds \eqref{eq:l2} then follow from the adjoint of Proposition \ref{prop:l2} (note that $U^\pm(0)^* g$ is supported on the support of $a$ and so the factor $\frac{1}{2\pi i |\xi|}$ is harmless).  To prove the remaining estimates, it suffices to prove the bounds
\begin{align}
\| U^\pm(0) U^\pm(0)^* f - f \|_2 &\lesssim \eps^{\delta} \label{uu-1}\\
\| U^+(0) \frac{1}{2\pi i |\xi|} U^+(0)^* g - U^-(0) \frac{1}{2\pi i |\xi|} U^-(0)^* g \|_2 &\lesssim \eps^{\delta} \label{uu-2}\\
\| \dot U^\pm(0) \frac{1}{\pm 2\pi i |\xi|} U^\pm(0)^* g - g \|_2 &\lesssim \eps^{\delta} \label{uu-3}\\
\| \dot U^+(0) U^+(0)^* f + \dot U^-(0) U^-(0)^* f \|_2 &\lesssim \eps^{\delta}.\label{uu-4}
\end{align}
for both choices of sign $\pm$.

Consider \eqref{uu-1}.  From \eqref{eq:uu-form} we have
$$
U^\pm(0) U^\pm(0)^* f(x) := \int (\int e^{2\pi i (\Psi_\pm(0,x,\xi) - \Psi_\pm(0,y,\xi))} e^{2\pi i (x-y) \cdot \xi} a^2(\xi)\ d\xi) f(y)\ dy.
$$
On the other hand, from the Fourier inversion formula, the Fourier support property of $f$, and the definition of $a$, we have
$$
f(x) := \int (\int e^{2\pi i (x-y) \cdot \xi} a^2(\xi)\ d\xi) f(y)\ dy.
$$
It thus suffices to prove the estimate
$$
\|\int (\int (e^{2\pi i (\Psi_\pm(0,x,\xi) - \Psi_\pm(0,y,\xi))}-1) e^{2\pi i (x-y) \cdot \xi} a^2(\xi)\ d\xi) f(y)\ dy\|_2 \lesssim \eps^{\delta}.
$$
We will in fact prove the slightly more general estimate
\begin{equation}\label{eq:osc}
\|\int (\int (e^{2\pi i (\Psi_\pm(0,x,\xi) - \Psi_\pm(0,y,\xi))}-1) e^{2\pi i (x-y) \cdot \xi} \eta(\xi)\ d\xi) F(y)\ dy\|_2 \lesssim \eps^{\delta} \|F\|_2
\end{equation}
for any bump function $\eta$ supported on the annulus $\{ 2^{-3} \leq |\xi| \leq 2^3 \}$ and any $F$.

By Schur's test and the self-adjointness of the kernel, it will suffice to prove the kernel bounds
\begin{align*}
|\int (e^{2\pi i (\Psi_\pm(0,x,\xi) - \Psi_\pm(0,y,\xi))}-1) e^{2\pi i (x-y) \cdot \xi} \eta(\xi)\ d\xi| \lesssim &K_0(x,y)+ \min( \eps |x-y|, \langle x-y\rangle^{-n-\sqrt{\delta}}),\\
\sup_x \int K_0(x,y)\, dy\les \epsilon
\end{align*}
since the right-hand side of the first line 
has an $x$ or $y$ integral of $O(\eps^\delta)$ (if $\delta$ is 
chosen sufficiently small), as can be seen by dividing into regions $|x-y| > \eps^{-\sqrt{\delta}}$ and $|x-y| \leq \eps^{-\sqrt{\delta}}$.

The bound of $\langle x-y \rangle^{-n-\sqrt{\delta}}$, as well as the bound for the 
kernel $K_0$, comes from \eqref{k-bound}, \eqref{k-bound-2} 
(note that the contribution of the $-1$ term can be dealt with directly by stationary phase).  To obtain the $\eps |x-y|$ term, simply observe from the mean-value theorem that
$$ |e^{2\pi i (\Psi_\pm(0,x,\xi) - \Psi_\pm(0,y,\xi))}-1| \lesssim |x-y| \| \nabla_x \Psi_\pm(0,x,\xi) \|_{L^\infty_x}.$$
But from \eqref{eq:qr} and Lemma \ref{lem:psid-1} we have
$$ \| \nabla_x \Psi_\pm \|_{L^\infty_t L^\infty_x} \lesssim \eps,$$
and the claim follows.

Now we prove \eqref{uu-2}.  A calculation similar to \eqref{eq:uu-form} shows that
$$
U^\pm(0) \frac{1}{2\pi i |\xi|} U^\pm(0)^* g(x) := \int (\int e^{2\pi i (\Psi_\pm(0,x,\xi) - \Psi_\pm(0,y,\xi))} e^{2\pi i (x-y) \cdot \xi} \frac{a^2(\xi)}{2\pi i |\xi|}\ d\xi) g(y)\ dy.
$$
Thus to prove \eqref{uu-2} it suffices by the triangle inequality to show that
$$
\| \int (\int (e^{2\pi i (\Psi_\pm(0,x,\xi) - \Psi_\pm(0,y,\xi))}-1) e^{2\pi i (x-y) \cdot \xi} \frac{a^2(\xi)}{2\pi i |\xi|}\ d\xi) g(y)\ dy \|_2 \lesssim \eps^{\delta}.$$
But this follows from \eqref{eq:osc} and \eqref{eq:fg-norm}.

The estimate \eqref{uu-3} follows from \eqref{uu-1} and
\eqref{eq:t-small} (as well as the adjoint of Proposition \ref{prop:l2}).  The estimate \eqref{uu-4} follows in a similar manner from \eqref{uu-2}.
\end{proof}

\section{Decay and Strichartz estimates for the parametrix}\label{sec:decay}

Fix $f,g$ as in Proposition \ref{Strichartz-0}, and let $h_+$, $h_-$ be the functions given by Proposition \ref{prop:unitary}.  We define the parametrix $\phi$ to be
$$ \phi(t) := \sum_\pm U_\pm(t) h_\pm.$$
By Proposition \ref{prop:unitary}, this parametrix obeys \eqref{eq:data-approx}.  In this section we now prove the Strichartz estimate \eqref{strichartz-covariant-localized}.  In light of Proposition \ref{prop:unitary}, it will suffice to prove the estimates
\begin{equation}\label{eq:fourier-strichartz} 
\| \nabla_{x,t} U_\pm(t) h\|_{L^\infty_t L^2_x} + \| \nabla_{x,t} U_\pm(t) h \|_{L^2_t L^{p_*}_x} \lesssim \|h\|_2
\end{equation}
on the spacetime slab $I \times \R^n$, for any function $h(\xi)$.  The former bound follows immediately from Proposition \ref{prop:l2-deriv}, so it suffices to prove the latter.

The first objective is to prove the (slightly simpler) Strichartz estimate without the $\nabla_{x,t}$ derivative:

\begin{proposition}\label{bs} We have
\begin{equation}\label{eq:baby-strichartz}
\|  U_\pm(t) h \|_{L^2_t L^{p_*}_x} \lesssim \|h\|_2.
\end{equation}
\end{proposition}

\begin{proof}
To obtain this Strichartz estimate we invoke the abstract Strichartz theorem from (\cite{tao:keel}, Theorem 1.2).  This theorem allows us to deduce \eqref{eq:baby-strichartz} from the energy estimate in Proposition \ref{prop:l2}, as long as we have the dispersive estimate
\begin{equation}\label{eq:dispersive}
 \| U_\pm(t) U_\pm^*(s) f \|_{L^\infty_x} \lesssim \langle t-s \rangle^{-\frac{n-1}{2} + \delta^2} \|f\|_{L^1_x}
\end{equation}
for all (Schwartz) $f$ and all $t,s \in \R$.

Fix $t,s$; we allow the functions we define to depend implicitly on $t$ and $s$.  
From \eqref{eq:ut-def} we see that
$$
U_\pm(t) U_\pm^*(s) f(x) = \int \int e^{2\pi i (\Psi_\pm(t,x,\xi) - \Psi_\pm(s,y,\xi))} e^{2\pi i (x-y) \cdot \xi} e^{\pm 2\pi i (t-s) |\xi|} a^2(\xi)\ d\xi f(y)\ dy$$
so it will suffice to prove that
\begin{equation}\label{eq:decay-est}
 |\int e^{2\pi i (\Psi_\pm(t,x,\xi) - \Psi_\pm(s,y,\xi))} e^{2\pi i (x-y) \cdot \xi} e^{\pm 2\pi i (t-s) |\xi|} a^2(\xi)\ d\xi | \lesssim 
\langle t-s \rangle^{-\frac{n-1}{2} + \delta^2}
\end{equation}
for all $x,y$.  We divide this into three cases: the \emph{strongly timelike case} $|x-y| \leq \frac{1}{2} |t-s|$,
the \emph{strongly spacelike case} $|x-y| \geq 2 |t-s|$, and the \emph{intermediate case} $\frac{1}{2}|t-s|
< |x-y| < 2|t-s|$.  Also we may assume that $|t-s| \gg 1$, since the claim follows trivially otherwise by
estimating everything by absolute values.

\divider{The strongly timelike case $|x-y| \leq \frac{1}{2}|t-s|$.}

This is the easiest case, as one can obtain very rapid decay purely from radial estimates.  Indeed, from 
integration by parts we see that
$$ |\int e^{2\pi i (x-y) \cdot \omega r} e^{\pm 2\pi i (t-s) r} a^2(r\omega)\ r^{n-1}\ dr | \lesssim |t-s|^{-100n}$$
for every $\omega \in S^{n-1}$, since the phase $(x-y) \cdot \omega \pm (t-s)$ has magnitude $\sim |t-s|$.  Integrating this in $\omega$ one obtains \eqref{eq:decay-est} (note that the phase $\Psi_\pm(t,x,\xi)$ depends only on the angular component $\omega = \xi/|\xi|$ of $\xi$ and not on the radial component).

\divider{The strongly spacelike case $|x-y| \geq 2|t-s|$.}

This case will resemble the argument used to prove \eqref{schur} (which corresponds to the special case
$t=s$ of the strongly spacelike case, although the bound required for \eqref{schur} was slightly different).  Fix $x,y,t,s$.  As in the proof of \eqref{schur}, we shall use
\eqref{eq:psi-dyadic} to  split the phase $\Psi_\pm(t,x,\xi) - \Psi_\pm(s,y,\xi)$ into two pieces,
\begin{equation}\label{eq:psi-decomp}
 \Psi_\pm(t,x,\xi) - \Psi_\pm(s,y,\xi) = \Phi_{smooth}(x,y,\xi) + \Phi_{small}(x,y,\xi),
\end{equation}
where 
\begin{equation}\label{eq:phi-small-def}
\Phi_{small}(x,y,\xi) := \Psi_\pm^{<\theta_*}(t,x,\xi) - \Psi_\pm^{<\theta_*}(s,y,\xi)
\end{equation}
and
\begin{equation}\label{eq:phi-smooth-def}
\Phi_{smooth}(x,y,\xi) := \Psi_\pm^{\geq \theta_*}(t,x,\xi) - \Psi_\pm^{\geq \theta_*}(s,y,\xi),
\end{equation}
$\theta_*$ is the angle $\theta_* := |x-y|^{-1+2\delta}$, and $\Psi_\pm^{<\theta_*}$ and $\Psi_\pm^{\geq \theta_*}$
were defined in \eqref{psi-less}, \eqref{psi-more}. To prove the claim \eqref{eq:decay-est}, we first use \eqref{eq:psi-decomp} to estimate
$$ e^{2\pi i (\Psi_\pm(t,x,\xi) - \Psi_\pm(s,y,\xi))} = e^{2\pi i (\Phi_{smooth}(\xi))} + O(|\Phi_{small}(\xi)|),$$
and then conclude that it will suffice to show the bounds
\begin{equation}\label{eq:small}
 \int |\Phi_{small}(x,y,\xi)| a^2(\xi)\ d\xi \lesssim |x-y|^{-\frac {n-1}2 }
\end{equation}
and
\begin{equation}\label{eq:pre-smooth}
 |\int e^{2\pi i \Phi_{smooth}(x,y,\xi)} e^{2\pi i (x-y) \cdot \xi} e^{\pm 2\pi i (t-s) |\xi|} a^2(\xi)\ d\xi | \lesssim 
\langle t-s \rangle^{-\frac{n-1}{2}}.
\end{equation}

\divider{Proof of \eqref{eq:small}.}

We first prove \eqref{eq:small}, which is proven similarly to \eqref{eq:Phi-small-bound}.  By the 
mean value theorem and the strongly spacelike hypothesis we have
$$ |\Phi_{small}(\xi)| \lesssim |x-y| \| \nabla_{x,t} \Psi_\pm^{<\theta_*}(t') \|_{L^\infty_x}$$
for some time $t'$ between $t$ and $s$.
By \eqref{psi-less} and the triangle inequality, the right-hand side is bounded by
$$ \lesssim |x-y| 
\sum_{k < -5} 
\sum_{j \geq 1: 2^{\sigma k + j} < \theta_*} 
\| \nabla^2_{x,t} (2^k  2^{\sigma k + j})^{-2} \tilde \Pi_{\omega, 2^{\sigma k + j}} \tilde P_k \underline{A}(t') \cdot\omega \|_{L^\infty_x}.$$
Thus by the Fourier inversion formula and Lemma \ref{lem:Coulomb}, we can bound the previous expression by
\begin{align*} & \lesssim |x-y| 
\sum_{k < -5} 
\sum_{j \geq 1: 2^{\sigma k + j} < \theta_*} (2^k 2^{\sigma k + j})^{-2}
\int_{|\zeta| \sim 2^k; \angle(\zeta,\omega) \sim 2^{\sigma k + j}} |\widehat{\nabla^2_{x,t} \underline{A}}(t',\zeta)\cdot \omega|\ d\zeta\\
& \lesssim |x-y| 
\sum_{k < -5} 
\sum_{j \geq 1: 2^{\sigma k + j} < \theta_*} 2^{-2k} 2^{-(\sigma k + j)}
\int_{|\zeta| \sim 2^k; \angle(\zeta,\omega) \sim 2^{\sigma k + j}} |\widehat{\nabla^2_{x,t} \underline{A}}(t',\zeta)|\ d\zeta
\end{align*}
Thus we can bound the left-hand side of \eqref{eq:small} by
$$ \lesssim |x-y| 
\sum_{k < -5} 
\sum_{j \geq 1: 2^{\sigma k + j} < \theta_*} 2^{-2k} 2^{-(\sigma k + j)}
\int_{S^{n-1}}
\int_{|\zeta| \sim 2^k; \angle(\zeta,\omega) \sim 2^{\sigma k + j}} |\widehat{\nabla^2_{x,t} \underline{A}}(t',\zeta)|\ d\zeta
d\omega.
$$
Interchanging the integrals, and then performing the $\omega$ integral, this is bounded by
$$ \lesssim |x-y| 
\sum_{k < -5} 
\sum_{j \geq 1: 2^{\sigma k + j} < \theta_*} 2^{-2k}(2^{\sigma k + j})^{n-2}
\int_{|\zeta| \sim 2^k} |\widehat{\nabla^2_{x,t} \underline{A}}(t',\zeta)|\ d\zeta.
$$
By Plancherel, Cauchy-Schwarz, and \eqref{eq:a-sob} we have that
$$ \int_{|\zeta| \sim 2^k} |\widehat{\nabla^2_{x,t} \underline{A}}(t',\zeta)|\ d\zeta \lesssim \eps 2^{3k} \lesssim 2^{3k}$$
and so we can bound the left-hand side of \eqref{eq:small} by
$$ \lesssim |x-y| 
\sum_{k < -5} 
\sum_{j \geq 1: 2^{\sigma k + j} < \theta_*} 2^{3k} 2^{-2k} (2^{\sigma k + j})^{n-2}.$$
Observe that the inner sum is vacuous unless $2^{\sigma k} < \theta_*$.  If we then perform the $j$ sum, we can estimate the previous by
$$ \lesssim |x-y| 
\sum_{k: 2^{\sigma k} < \theta_*}  2^k \theta_*^{n-2}.$$
Summing over all $k$ (using the hypothesis $n \geq \n$), we can bound the previous by
$$ \lesssim |x-y| \theta_*^{n-2 + 1/\sigma} = |x-y|^{1 - (n-2 + 1/\sigma)(1-2\delta)}.$$
This is acceptable if
$$ 1 - (n-2 + \frac{1}{\sigma}) < -n/2+1/2$$
and $\delta$ is sufficiently small.  But this inequality follows automatically if $n\ge 5$.  This proves
\eqref{eq:small}.

\divider{Proof of \eqref{eq:pre-smooth}.}

We now prove \eqref{eq:pre-smooth}.  As in the proof of \eqref{eq:Phi-smooth-bound}, it will suffice to
prove the smoothness estimate
\begin{equation}\label{smooth}
 |\nabla_\xi^K \Phi_{smooth}(x,y,\xi)| \lesssim |x-y|^{K(1-\delta)}
\end{equation}
for all $K \geq 1$ and all $\xi$ in the support of $a$.  Note that the phase factor $(t-s)|\xi|$ can be grouped
with the more rapidly oscillating phase $(x-y) \cdot \xi$ and will not significantly affect the stationary
phase analysis since we are in the strongly spacelike case $|x-y| \geq 2 |t-s|$.

Now we prove \eqref{smooth}, which is a similar estimate to \eqref{smooth-est}.  Fix $\xi$ in the support of $a$.  
As in the proof of \eqref{smooth-est} we observe that we have the two estimates
$$ |f(t,x) - f(s,y)| \lesssim \min_{m=0,1} |x-y|^m \| \nabla^m_{x',t'} f(t') \|_{L^\infty_{x'}}$$
for some time $t'$ between $t$ and $s$.  We can thus estimate the left-hand side of \eqref{smooth} as
$$
\lesssim \sum_{k < -5} \sum_{j \geq 1: 2^{\sigma k + j} \geq \theta_*} 
\min_{m=0,1} \sup_{t'} |x-y|^m \| \nabla_\xi^K \nabla^m_{x',t'} (c_\pm(\omega) \cdot \nabla_{x',t'} (2^k 2^{\sigma k + j})^{-2} \tilde \Pi_{\omega, 2^{\sigma k + j}} \tilde P_k \underline{A}(t')) \cdot\omega\|_{L^\infty_{x'}}.$$
Thus, by arguing as in the proof of \eqref{smooth-est}, we can bound this expression by
$$
\lesssim \sum_{k < -5} \sum_{j \geq 1: 2^{\sigma k + j} \geq \theta_*} 
\min_{m=0,1} \sup_{t'} \int_{|\zeta| \sim 2^k; \angle(\zeta,\omega) \sim 2^{\sigma k + j}}
(2^{-\sigma k - j})^K |x-y|^m 2^{-2k} 2^{-(\sigma k + j)}
|\widehat{\nabla^{m+1}_{x',t'} \underline{A}(t')}(\zeta)|\ d\zeta.$$
We may of course restrict the $j$ summation to the region where $2^{\sigma k + j} \lesssim 1$, since the integral is vacuous otherwise.  But then by \eqref{repeat}, with the additional
help of  \eqref{eq:a-sob}, as before we may estimate the above expression by \eqref{repeat-2},
and so \eqref{smooth} then follows by repeating the proof of \eqref{smooth-est}.
This concludes the treatment of the strongly spacelike case.

\divider{The intermediate case $\frac{1}{2}|t-s| \leq |x-y| \leq 2|t-s|$.}

It remains to consider the intermediate regime $|x-y| \sim |t-s|$, which is neither strongly timelike nor strongly spacelike.

Without loss of generality we may write $x-y = |x-y| e_1$, where $e_1$ is the first standard basis vector.  We then partition the symbol $a^2(\xi)$ as
$$ a^2(\xi) = b_0(\xi) + \sum_{j: 1 < 2^m \lesssim |t-s|^{1/2}} b_m(\xi)$$
where $b_0(\xi)$ is a non-negative bump function adapted to the sector $\{ |\xi| \sim 1, \angle(\xi, e_1) \lesssim |t-s|^{-1/2} \}$, and $b_m(\xi)$ is a non-negative bump function adapted to the sector $\{ |\xi| \sim 1, \angle(\xi,e_1) \sim 2^m |t-s|^{-1/2} \}$.  The contribution of the $b_0(\xi)$ term is certainly $O(|t-s|^{-(n-1)/2})$, just by replacing everything by absolute values, so it will suffice to show that
\begin{equation}\label{eq:jsum}
|\int e^{2\pi i (\Psi_\pm(t,x,\xi) - \Psi_\pm(s,y,\xi))} e^{2\pi i (x-y) \cdot \xi} e^{\pm 2\pi i (t-s) |\xi|} b_m(\xi)\ d\xi | \lesssim 2^{-\delta m} |t-s|^{-(n-1)/2 + \delta^2}
\end{equation}
for each $1 < 2^m \lesssim |t-s|^{1/2}$.

Fix $m$.  We repeat our strategy used to treat the strongly spacelike case.  Namely, we again perform the decomposition \eqref{eq:psi-decomp}, where $\Phi_{small}$ and $\Phi_{smooth}$ are again defined by \eqref{eq:phi-small-def} and \eqref{eq:phi-smooth-def} respectively.  This time, however, the critical angle $\theta_*$ will be selected according to the formula
$$ \theta_* := 2^{-m} |t-s|^{-\frac{1}{2} + \delta^3}.$$

We shall prove the bounds
\begin{equation}\label{eq:small-j}
 \int |\Phi_{small}(\xi)| b_m(\xi)\ d\xi \lesssim 2^{-\delta m} |t-s|^{-(n-1)/2 + \delta^2}
\end{equation}
and 
\begin{equation}\label{eq:smooth-j}
 |\nabla_\xi^K \Phi_{smooth}(\xi)| \lesssim (\log |t-s|/\theta_*)^K
\end{equation}
for all integers $K \geq 1$ (cf. \eqref{eq:small}, \eqref{smooth}).  If we assume these bounds, then by arguing as in the strongly spacelike case, we may use \eqref{eq:small-j} to replace the left-hand side of \eqref{eq:jsum} by
\begin{equation}\label{eq:jsum-smooth}
 |\int e^{2\pi i \Phi_{smooth}(\xi)} e^{2\pi i (x-y) \cdot \xi} e^{\pm 2\pi i (t-s) |\xi|} b_m(\xi)\ d\xi|.
\end{equation}
Then observe from \eqref{eq:smooth-j} that
$$ | \nabla_\xi^K e^{2\pi i \Phi_{smooth}(\xi)}| \lesssim (\log |t-s|/\theta_*)^K $$
for all $K \geq 0$, while from the definition of $b_m$ we have
$$ | \nabla_\xi^K b_m(\xi)| \lesssim (2^m |t-s|^{-1/2})^{-K}  $$.
Thus by the Leibnitz rule
$$ | \nabla_\xi^K (e^{2\pi i \Phi_{smooth}(\xi)} b_m(\xi))| \lesssim (\log |t-s|/\theta_* + 2^{-m} |t-s|^{1/2})^K.$$
On the other hand, observe that
$$ |\nabla_\xi( (x-y) \cdot \xi \pm (t-s) |\xi| )| = |(x-y) \pm (t-s) \frac{\xi}{|\xi|}| \gtrsim |t-s| |\xi - \xi_1 e_1| \gtrsim |t-s|^{1/2} 2^m.$$
Thus by repeated integration by parts, we can estimate \eqref{eq:jsum-smooth} by
$$ \lesssim (|t-s|^{1/2} 2^m)^{-K} (\log |t-s|/\theta_*+ 2^{-m} |t-s|^{1/2})^K = |t-s|^{-\delta^3 K} (\log |t-s|)^K + 2^{-2mK}$$
for any integer $K \geq 0$.  Setting $K$ sufficiently large and using the size of the support 
of $b_m(\xi)$, we obtain that
$$
\eqref{eq:jsum-smooth}\les (|t-s|^{-\delta^3 K} (\log |t-s|)^K + 2^{-2mK}) 2^{m(n-1)}
|t-s|^{-(n-)/2}\les 2^{-\delta m} |t-s|^{-(n-1)/2+\delta^2}
$$
as desired.
It remains to prove \eqref{eq:small-j} and \eqref{eq:smooth-j}.  The estimate \eqref{eq:smooth-j} follows from \eqref{smooth} (observe that the proof of \eqref{smooth} did not use anything about the choice of $\theta_*$), so we only need to show \eqref{eq:small-j}.  Arguing exactly as in the proof of \eqref{eq:small}, we have
$$ |\Phi_{small}(\xi)| \lesssim |x-y| 
\sum_{k < -5} 
\sum_{j \geq 1: 2^{\sigma k + j} < \theta_*} 2^{-2k} 2^{-(\sigma k + j)}
\int_{|\zeta| \sim 2^k; \angle(\zeta,\omega) \sim 2^{\sigma k + j}} |\widehat{\nabla^2_{x,t} \underline{A}}(t,\zeta)|\ d\zeta.
$$
By hypothesis we have $|x-y| \sim |t-s|$.
Thus we can bound the left-hand side of \eqref{eq:small} by
$$ \lesssim |t-s| 
\sum_{k < -5} 
\sum_{j \geq 1: 2^{\sigma k + j} < \theta_*} 2^{-2k} 2^{-(\sigma k + j)}
\int_{\angle(\omega,e_1) \lesssim 2^m |t-s|^{-1/2}}
\int_{|\zeta| \sim 2^k; \angle(\zeta,\omega) \sim 2^{\sigma k + j}} |\widehat{\nabla^2_{x,t} \underline{A}}(t,\zeta)|\ d\zeta
d\omega.
$$
Observe that $2^{\sigma k + j} < \theta_*$ and $\theta_* + 2^m |t-s|^{-1/2} \lesssim 2^m |t-s|^{-1/2+\delta^3}$.  Thus the $\zeta$ variable is constrained to the sector where $\angle(\zeta,e_1) \lesssim 2^m |t-s|^{-1/2+\delta^3}$.  Interchanging the integrals, and then performing the $\omega$ integral, this is bounded by
$$ \lesssim |t-s| 
\sum_{k < -5} 
\sum_{j \geq 1: 2^{\sigma k + j} < \theta_*} 2^{-2k} 2^{-(\sigma k + j)} 
(2^{\sigma k + j})^{n-1}
\int_{|\zeta| \sim 2^k; \angle(\zeta,e_1) \lesssim 2^m |t-s|^{-1/2+\delta^3}} |\widehat{\nabla^2_{x,t} \underline{A}}(t,\zeta)|\ d\zeta.
$$
By Plancherel, Cauchy-Schwarz, and \eqref{eq:a-sob} we have that
$$ \int_{|\zeta| \sim 2^k; \angle(\zeta,e_1) \lesssim \theta_*} |\widehat{\nabla^2_{x,t} \underline{A}}(t,\zeta)|\ d\zeta \lesssim 2^{3k} (2^m |t-s|^{-1/2+\delta^3})^{(n-1)/2},$$
and so we can bound the left-hand side of \eqref{eq:small} by
$$ \lesssim |t-s| (2^m |t-s|^{-1/2+\delta^3})^{(n-1)/2}
\sum_{k < -5} 
\sum_{j \geq 1: 2^{\sigma k + j} < \theta_*} 2^{k} (2^{\sigma k + j})^{n-2}.$$
Observe that the inner sum is vacuous unless $2^{\sigma k} < \theta_*$.  If we then perform the $j$ sum, we can estimate the previous by
$$ \lesssim |t-s| (2^m |t-s|^{-1/2+\delta^3})^{(n-1)/2}
\sum_{k: 2^{\sigma k} < \theta_*}  2^k \theta_*^{n-2}.$$
Summing over all $k$ (using the hypothesis $n \geq \n$), we can bound the previous by
$$ \lesssim |t-s| (2^m |t-s|^{-1/2+\delta^3})^{(n-1)/2} \theta_*^{n-2 + 1/\sigma}
= |t-s|^{1 - (\frac{3n-5}{2} + 1/\sigma)(\frac{1}{2}-\delta^3)} 
2^{(\frac{3-n}{2}-1/\sigma) m}.$$
This is acceptable if
$$ 1 - \frac{1}{2}(\frac{3n-5}{2} + \frac{1}{\sigma}) < -\frac{n-1}{2}$$
and $\delta$ is sufficiently small.  This inequality is implied automatically for $n\ge 7$, while for $n=\n$ it requires that $\si <2$. This  follows from \eqref{sigma-bounds}.
This proves \eqref{eq:small-j}.

This completes the proof of the dispersive estimate \eqref{eq:dispersive}, and hence \eqref{eq:baby-strichartz}.  
\end{proof}

We now re-introduce the spacetime derivative, and prove the estimate
\begin{equation}\label{eq:adult-strichartz}
\|  \nabla_{x,t} U_\pm(t) h \|_{L^2_t L^{p_*}_x} \lesssim \|h\|_2.
\end{equation}

This will be achieved by an argument similar to Proposition \ref{prop:l2-deriv}.
From \eqref{eq:ut-def} we have
$$ \partial_t U_\pm(t) := \int e^{2\pi i\Psi_\pm(t,x,\xi)} e^{2\pi i x \cdot \xi} e^{\pm 2\pi i t|\xi|} (2\pi i \partial_t \Psi_\pm \pm 2\pi i |\xi|) h_\pm(\xi) a(\xi)\ d\xi $$
and
$$ \nabla_x U_\pm(t) := \int e^{2\pi i\Psi_\pm(t,x,\xi)} e^{2\pi i x \cdot \xi} e^{\pm 2\pi i t|\xi|} (2\pi i \nabla_{x} \Psi_\pm \pm 2\pi i \xi) h_\pm(\xi) a(\xi)\ d\xi.$$
But by Lemma \ref{lem:psid-1} and \eqref{eq:factor}, we have
$$ \| 2\pi i \partial_t \Psi_\pm \pm 2\pi i |\xi| \|_{D(L^\infty_t L^\infty_x)},
\| 2\pi i \nabla_{x} \Psi_\pm \pm 2\pi i \xi \|_{D(L^\infty_t L^\infty_x)} \lesssim 1,$$
and so \eqref{eq:adult-strichartz} follows from \eqref{eq:baby-strichartz}.  As discussed at the beginning of this section, this proves \eqref{eq:fourier-strichartz} and hence \eqref{strichartz-covariant-localized}.

\section{The parametrix is accurate}\label{sec:accuracy}

In the last few sections we have constructed a parametrix for the modified covariant wave equation $\Box'_{\underline{A}} \phi = 0$, shown that it approximates the initial data at time 0, and obeys the desired Strichartz estimates.  To complete the proof of Proposition \ref{localized-Strichartz}, all we have to do now is show that the parametrix is indeed an approximate solution to the modified covariant wave equation, in the sense of \eqref{eq:forcing-approx}.  In light of \eqref{eq:ut-box} and Proposition \ref{prop:unitary}, it suffices to show that
\begin{equation}\label{eq:omega-coupled}
 \| \int \Omega_\pm(t,x,\xi)
e^{2\pi i\Psi_\pm(t,x,\xi)} e^{2\pi i x \cdot \xi} e^{\pm 2\pi i t|\xi|} h(\xi) a(\xi)\ d\xi \|_{L^1_t L^2_x} \lesssim \eps \| h\|_2
\end{equation}
for any (Schwartz) function $h(\xi)$ and either choice of sign $\pm$.

We use the splitting in Lemma \ref{lem:psid-2}.  Clearly it suffices to prove
$$ \| \int \Omega^1_\pm(t,x,\xi)
e^{2\pi i\Psi_\pm(t,x,\xi)} e^{2\pi i x \cdot \xi} e^{\pm 2\pi i t|\xi|} h(\xi) a(\xi)\ d\xi \|_{L^1_t L^2_x} \lesssim \eps \| h\|_2$$
and
$$ \| \int \Omega^2_\pm(t,x,\xi)
e^{2\pi i\Psi_\pm(t,x,\xi)} e^{2\pi i x \cdot \xi} e^{\pm 2\pi i t|\xi|} h(\xi) a(\xi)\ d\xi \|_{L^1_t L^2_x} \lesssim \eps \| h\|_2.$$
But the first estimate follows from the bound $\| \Omega^1_\pm\|_{D(L^1_t L^\infty_x)} \lesssim \eps$ and Proposition \ref{prop:l2}, while the second estimate follows from the bound $\| \Omega^2_\pm \|_{D(L^2_t L^q_x)} \lesssim \eps$ and Proposition \ref{bs}.

This completes the proof of Proposition \ref{localized-Strichartz}, and hence Theorem \ref{main}.


\begin{thebibliography}{10}

\bibitem{bc}
H. Bahouri, J-Y. Chemin, 
\emph{\'Equations d'ondes quasilin\'eaires et les inegalites de Strichartz}, Amer. J. Math.
\textbf{121} (1999), 1337--1377. 

\bibitem{cuccagna} 
S. Cuccagna, \emph{On the local existence for the Maxwell-Klein-Gordon system in $\R^{3+1}$}, 
Comm. Partial Differential Equations. \textbf{24} (1999), 851--867. 

\bibitem{em}
D. Eardley, V. Moncrief, \emph{The global existence of Yang-Mills-Higgs fields in $\R^{3+1}$}, 
Comm. Math. Phys. \textbf{83} (1982), 171--212. 

\bibitem{tao:keel}
M. Keel, T. Tao, \emph{Endpoint Strichartz Estimates}, Amer. Math. J. \textbf{120} (1998), 955--980.

\bibitem{keel:mkg}
M. Keel, T. Tao, \emph{Global well-posedness for large data for the Maxwell-Klein-Gordon equation below the energy norm}, in preparation. 

\bibitem{kl-mac}
S. Klainerman, M. Machedon, \emph{On the Maxwell-Klein-Gordon equation with finite energy}, Duke Math. J. \textbf{74} (1994), 19--44. 

\bibitem{kr:wavemap}
S. Klainerman, I. Rodnianski, \emph{On the global regularity of wave maps in the critical Sobolev norm}, IMRN \textbf{13} (2001), 656--677.

\bibitem{kr:einstein}
S. Klainerman, I. Rodnianski, \emph{Rough solution for the Einstein Vacuum equations}, preprint.

\bibitem{kl-tar}
S. Klainerman, D. Tataru, \emph{On the optimal regularity for Yang-Mills equations in $\R^{4+1}$}, J. Amer. Math. Soc. \textbf{12} (1999), 93--116. 

\bibitem{koch}
H. Koch, N. Tzvetkov, \emph{On the local well-posedness of the Benjamin-Ono equation in $H^s(\R)$}, IMRN  2003:26. 

\bibitem{mac}
M. Machedon, J. Sterbenz, \emph{Optimal well-posedness for the Maxwell-Klein-Gordon equations in 3+1 dimensions}, preprint. 

\bibitem{nahmod}
A. Nahmod, A. Stefanov, K. Uhlenbeck, \emph{On the well-posedness of the wave map problem in high dimensions}, Comm. Anal. Geom. \textbf{11} (2003), 49--83. 

\bibitem{selberg:mkg}
S. Selberg, \emph{Almost optimal local well-posedness of the Maxwell-Klein-Gordon equations on $\R^{1+4}$}, preprint.

\bibitem{shsw:wavemap}
J. Shatah, M. Struwe, \emph{The Cauchy problem for wave maps}, IMRN \textbf{11} (2002) 555--571.

\bibitem{smith}
H. Smith, D. Tataru, \emph{Sharp local well-posedness results for the nonlinear wave equation}, preprint. 

\bibitem{sterbenz}
J. Sterbenz, \emph{Null concentration, scale invariance, and global regularity for quadratic non-linear wave equations I.  Critical Besov spaces in high dimensions}, preprint.

\bibitem{tao:wavemap1}
T. Tao, \emph{Global regularity of wave maps I.  Small critical Sobolev norm in high dimension}, IMRN \textbf{7} (2001), 299--328.

\bibitem{tataru}
D. Tataru, \emph{Rough solutions for the wave maps equation}, preprint. 

\end{thebibliography}
\end{document}